\documentclass[12pt]{amsart}
\usepackage{amssymb}
\usepackage{amsfonts}
\usepackage{amscd}
\usepackage{latexsym}

\newcounter{cs}
\stepcounter{cs}
\newcounter{ds}
\stepcounter{ds}

\newcommand{\casos}{\begin{itemize}}
\newcommand{\fcasos}{\end{itemize}\setcounter{cs}{1}}

\newcommand{\End}{\mbox{\rm End}}
\newcommand{\Hom}{\mbox{\rm Hom}}
\newcommand{\soc}{\mbox{\rm Soc}}

\newcommand{\la}{\langle}
\newcommand{\ra}{\rangle}

\newcommand{\ol}{\overline}
\newcommand{\krat}{k_{\mbox{\rm rat}}\langle X\rangle}

\newcommand{\im}{\mbox{Im}}

\newtheorem{lem}{Lemma}[section]
\newtheorem{corol}[lem]{Corollary}
\newtheorem{theor}[lem]{Theorem}
\newtheorem{prop}[lem]{Proposition}

\theoremstyle{definition}

\newtheorem{defis}[lem]{Definitions}
\newtheorem{exem}[lem]{Example}
\newtheorem{exems}[lem]{Examples}
\newtheorem{nota}[lem]{Notation}
\newtheorem{rema}[lem]{Remark}
\newtheorem{remas}[lem]{Remarks}

\newcommand{\N}{\mathbb{N}}

\newcommand{\Z}{\mathbb{Z}}
\newcommand{\etil}{\tilde{e}}
\newcommand{\Etil}{\widetilde{E}}

\begin{document}

\title{$K_0$ of purely infinite simple regular rings}

\author{P. Ara}\address{Departament de Matem\`atiques, Universitat Aut\`onoma de
Barcelona, 08193, Bellaterra (Barcelona),
Spain}\email{para@mat.uab.es}
\author{K. R. Goodearl}\address{Department of Mathematics, University of California,
Santa Barbara, CA 93106, USA}\email{goodearl@math.ucsb.edu}
\author{E. Pardo}\address{Departamento de Matem\'aticas, Universidad de C\'adiz, Aptdo. 40,
11510 Puerto Real (C\'adiz), Spain}\email{enrique.pardo@uca.es}

\thanks{The first and third authors were partially supported by
MEC-DGESIC grant PB98-0873, and by the Comissionat per Universitats i
Recerca de la Generalitat de Catalunya. The second author was
partially supported by NSF grant DMS-9970159, and the third author by
PAI III grant FQM-298}

\begin{abstract}
We extend the notion of a purely infinite simple C*-algebra to the
context of unital rings, and we study its basic properties, specially
those related to K-Theory. For instance, if $R$ is a purely infinite
simple ring, then $K_0(R)^+= K_0(R)$, the monoid of isomorphism classes
of finitely generated projective $R$-modules is isomorphic to the
monoid obtained from $K_0(R)$ by adjoining a new zero element, and
$K_1(R)$ is the abelianization of the group of units of $R$. We develop
techniques of construction, obtaining new examples in this class in
the case of von Neumann regular rings, and we compute the Grothendieck
groups of these examples. In particular, we prove that every countable
abelian group is isomorphic to $K_0$ of some purely infinite simple
regular ring. Finally, some known examples are analyzed within this
framework.
\end{abstract}

\maketitle

\section*{Introduction}

In 1981, Cuntz \cite{C2} introduced the concept of a purely infinite
simple C*-algebra. This notion has played a central role in the
development of the theory of C*-algebras in the last two decades.
A large series of contributions, due to Blackadar, Brown, Lin,
Pedersen, Phillips, R\o rdam and Zhang, among others, reflect the
interest in the structure of such algebras. One of
the most important advances in the program of classifying separable
C*-algebras through K-Theory, proposed by Elliott in the early
seventies, was obtained in this context by Kirchberg \cite{Kirch} and
Phillips \cite{NCP}, who showed that nuclear separable unital purely
infinite simple C*-algebras are classified by K-theoretic
invariants.

In this work, we introduce a suitable definition of a purely infinite
simple ring. This notion agrees with that of Cuntz in the case of
C*-algebras. Moreover, a number of basic results, known in the
case of C*-algebras, also hold in the purely algebraic context. In
particular the algebraic $K_0$ and $K_1$ groups of a purely infinite
simple ring follow the same patterns as the corresponding topological
$K$ groups of purely infinite simple C*-algebras, found by Cuntz
in \cite{C2}.

Our goal is to use these concepts and ideas in order to advance
the knowledge of the K-theory of (von Neumann) regular rings. In
particular, we construct examples of purely infinite simple regular
rings whose
$K_0$ groups are cyclic of arbitrary order, and -- in the case when
$K_0$ is finite cyclic -- whose finitely generated projective modules
are free. By using these constructions, we build examples of purely
infinite simple regular rings whose
$K_0$ is any countable abelian group.

A key point of our work lies in the development of new techniques for
constructing examples. Our algebras are extensions of the algebras
$V_{1,n}$, first considered by Leavitt in \cite{Leavitt}. Recall that
the Leavitt algebra $V_{1,n}$ over a field $k$ is the $k$-algebra
having a universal right $V_{1,n}$-module isomorphism
$V^n_{1,n}\rightarrow V_{1,n}$. Our method for getting the extensions
of the $V_{1,n}$ works as follows. We consider a local subalgebra of
the algebra of noncommutative power series $k\la\la X\ra\ra$, where
$X=\{x_0,\dots ,x_n\}$, closed under a certain set of skew
derivations, and containing the free algebra $k\la X\ra$. Then we
construct our extensions of $V_{1,n}$ by universally inverting the
right $R$-module homomorphism $R\rightarrow R^{n+1}$ given by left
multiplication by the column $(x_0,\dots ,x_n)^T$. We prove that
$K_0$ of these algebras is always the cyclic group of order $n$.
Taking $R$ to be the full algebra of noncommutative power series one
obtains what could be thought of as a completion of $V_{1,n+1}$.
Taking the minimal choice for $R$, which turns to be the so-called
algebra of rational series, one obtains a kind of ``algebra of
fractions" of $V_{1,n+1}$. The latter algebra is shown to coincide
with a construction of Schofield \cite{SchDicks}, and also with a
construction of Rosenmann and Rosset \cite{R-R}. In particular it
enjoys various additional universal properties.

\medskip

We briefly outline the contents of the paper. In Section 1 we define
purely infinite simple rings, and we prove some of their elementary
properties. In Section 2 we give the basic patterns for the algebraic
$K_0$ and $K_1$ of purely infinite simple rings, in analogy with
Cuntz's results on the topological $K_0$ and $K_1$ of purely infinite
simple C*-algebras. Section 3 is devoted to introducing the
construction of skew polynomial rings with freely independent
indeterminates, which will be a fundamental technical device for the
rest of the paper. In Section 4 we outline some known results on
Leavitt's algebras, which can be considered as the earliest universal
examples in this class, and also as the algebraic relatives of
Cuntz's algebras. (The Cuntz algebra $\mathcal O _n$ is the
C*-completion of the Leavitt algebra $V_{1,n}$ over the field of
complex numbers.) Section 5 is the core of the paper: we develop here
the construction of our examples, we establish their structure, and we
compute their $K_0$ groups. In Section 6 we present the examples of
Schofield \cite{SchDicks}, as well as those of Rosenmann and Rosset
\cite{R-R}, and we show that they are isomorphic. Section \ref{isoQT} is
devoted to show that the examples in Section 6 are isomorphic to
particular cases of the examples presented in Section 5. As
a consequence, we see that these algebras can be obtained in two
different ways as universal localizations of free algebras.
Finally, we prove in Section 8 that any countable abelian group can be
realized as
$K_0$ of a purely infinite simple regular ring.

\medskip

Aside from a few noted exceptions, all rings and modules in this
paper will be assumed to be unital. In fact, most of our rings will
be algebras over a base field that we denote $k$. Ring and algebra
homomorphisms, except for embeddings of ideals, will also be assumed
to be unital. We will often write
$nA$ for the direct sum of $n$ copies of a module $A$, although in the
case of a ring $R$ we prefer the notation $R^n$ for the free right
$R$-module of rank $n$, whose elements we think of as column vectors.
The notations $k\la X\ra$ and $k\la\la X\ra\ra$, respectively, will
stand for the free algebra and the noncommutative formal power series
algebra over $k$ on a set $X$.

\section{Purely infinite simple rings}

We introduce the concept of a purely infinite simple ring, and
sketch some basic results on this topic. These results are analogous
to those obtained by Cuntz.

First, recall that if $R$ is a ring and $e,f\in R$ are idempotents,
we say that $e$ and $f$ are (Murray--von Neumann) {\it equivalent\/}
(denoted $e\sim f$) provided that there exist elements $x\in eRf$
and
$y\in fRe$ such that $xy=e$ and $yx=f$. This is equivalent to
demanding that $eR\cong fR$ as right $R$-modules. Recall also that
$e$ and $f$ are said to be {\it orthogonal\/} (denoted $e \perp f$)
provided that $ef=fe=0$. In that case, $e+f$ is an idempotent, and
$(e+f)R= eR \oplus fR$.

The following useful means of producing orthogonal decompositions of
idempotents is old and well known; we sketch the easy proof for
convenience.

\begin{lem}\label{pinflem1} Let $e$ be an idempotent in a ring $R$.
If there exist right ideals $A_i \subseteq R$ such that $eR=
A_1\oplus \cdots\oplus A_n$, then there exist pairwise orthogonal
idempotents $e_i\in R$ such that $e= e_1+ \dots+ e_n$ and $A_i=
e_iR$ for all $i$.
\end{lem}

\begin{proof} Due to the given direct sum decomposition, the
endomorphism ring of the module $eR$ contains orthogonal idempotents
$f_1,\dots,f_n$ such that $f_1+ \dots+ f_n$ equals the identity map
on $eR$ and each $f_i(eR)= A_i$. The desired idempotents $e_i$ are
the images of the $f_i$ under the canonical isomorphism
$\text{End}_R(eR) \rightarrow eRe$.
\end{proof}

\begin{defis} An idempotent $e$ in a ring $R$ is {\it infinite\/}
if there exist orthogonal idempotents $f,g\in R$ such that $e= f+g$
while $e\sim f$ and $g\ne 0$. In view of Lemma \ref{pinflem1}, $e$ is
infinite if and only if $eR$ is isomorphic to a proper direct summand
of itself, that is, $eR$ is a {\it directly infinite\/} module.

A  simple ring $R$ is said to be {\it purely infinite\/} if every
nonzero right ideal of $R$ contains an infinite idempotent. It will
follow from Theorem \ref{pinfthm} that this concept is left-right
symmetric.
\end{defis}

\begin{exems}\label{pinfexamples} The class of purely infinite
simple rings is rather large; we indicate some subclasses here:
\begin{enumerate}
\item[(a)] Many purely infinite simple C*-algebras are known; for
instance, see \cite{CA-D}, \cite{C1}, \cite{C2}, \cite{Kirch},
\cite{K-Ph}, \cite{KPRR},  \cite{LZ}, \cite{NCP},
\cite{rordamclassif}, \cite{Zhang}.
\item[(b)] If $V$ is an infinite dimensional vector space over $k$,
then $\End_k(V)$ modulo its unique maximal ideal $M$ is purely
infinite. This follows from the fact that if $f\in \End_k(V)
\setminus M$, then $\dim_k f(V)= \dim_k V$.
\item[(c)] More
generally, if $R$ is a regular, right self-injective ring without
nonzero directly finite central idempotents, and $M$ is any maximal
ideal of $R$, then $R/M$ is purely infinite. To see this, consider
any element $x\in R \setminus M$. Since $R/M$ is simple,
$\sum_{i=1}^n a_ixb_i -1 \in M$ for some $a_i,b_i \in R$, from which
we see that $R_R$ embeds in $n(xR)\oplus zR$ for some $z\in M$. By
general comparability \cite[Corollary 9.15]{vnrr}, there is a central
idempotent $e\in R$ such that $ezR \lesssim exR$ and $(1-e)xR
\lesssim (1-e)zR$; in particular, $(1-e)x \in M$. Since $x\notin M$, we
cannot have
$e\in M$, and so
$1-e\in M$. Now $eR$ is isomorphic to $exR$. For, since $e$ is a
directly infinite central idempotent, we have $(n+1)(eR)\cong eR$ by
\cite[Theorem 10.16]{vnrr}. Thus, as $eR\lesssim (n+1)(exR)\lesssim
(n+1)(eR)$, we have $(n+1)(eR)\cong (n+1)(exR)$ by \cite[Theorem
10.14]{vnrr}. Hence,
$eR\cong exR$ by \cite[Theorem 8.16(b)]{vnrr}. Since the class of $e$
in $R/M$ is $1$, we conclude that $R/M$ is isomorphic to $x(R/M)$.
\item[(d)] If $R$ is a directly infinite regular ring
and all (finitely generated) projective right $R$-modules are free,
then $R$ is a purely infinite simple ring. First, if $x$ is a
nonzero element of $R$, then the projective module $xR$ is free,
whence $R_R$ is isomorphic to a direct summand of $xR$, and
consequently $RxR=R$. Thus $R$ is simple. Further, since $R_R$ is
directly infinite, so is $xR$. Hence, all nonzero idempotents in
$R$ are infinite, and therefore $R$ is purely infinite.
\item[(e)] If $A$ is any directly infinite simple ring, then the
direct limit $B= \varinjlim M_{2^n}(A)$ (with block diagonal
transition maps) is purely infinite. To show this, let $z$ be a
nonzero element of $B$; then $z$ is the image of a nonzero element
$x\in M_{2^n}(A)$ for some $n$, and we can assume without loss of
generality that $n=0$. For a suitable $m$ we have that
$\sum_{i=1}^{2^m} a_ixb_i =1$ for some $a_i$, $b_i$ in $A$. If we use
the $a_i$ (respectively, the $b_i$) as the first row (respectively,
column) of $2^{m}\times 2^{m}$ matrices whose other entries are zero,
we obtain $a,b \in M_{2^{m}}(A)$ such that $ayb= e_{11}$, where $y$
is the image of $x$ in $M_{2^{m}}(A)$ and $e_{11}$ is the matrix
unit in the top left corner of $M_{2^{m}}(A)$. Since the idempotent
$1\in A$ is infinite, so is the idempotent $ayb$ in $M_{2^{m}}(A)$.
But $yba$ is an idempotent equivalent to $ayb$, so it too is
infinite. Its image in $B$ is an infinite idempotent contained in the
right ideal $zB$.
\end{enumerate}
\end{exems}

\begin{lem}\label{pinflem2} Let $R$ be a  simple ring, and let $P$
and $Q$ be finitely generated projective right $R$-modules. If $P$
is directly infinite, then there exists a nonzero right $R$-module
$A$ such that $P\cong Q\oplus A$.
\end{lem}

\begin{proof} By hypothesis, there exists a nonzero right
$R$-module $B$ such that $P\cong P\oplus B$, whence $P\cong P\oplus
mB$ for all $m\in\N $. Since $R$ is simple and $B$ is projective, $B$
is a generator in $\text{Mod-}R$; in particular, there exist $n\in\N$
and a right $R$-module $C$ such that $nB\cong Q\oplus C$. Therefore
$$P\cong P\oplus nB\cong Q\oplus (P\oplus C),$$ and $P\oplus C$ is
nonzero because $P\ne 0$.
\end{proof}

\begin{prop}\label{pinfprop} Suppose that $R$ is a  purely infinite
simple ring. Then all nonzero finitely generated projective right
$R$-modules are directly infinite; equivalently, all nonzero
idempotents in the matrix rings $M_n(R)$ are infinite.

Consequently, if $P$ and $Q$ are any nonzero finitely generated
projective right $R$-modules, then there exists a nonzero right
$R$-module $A$ such that $P\cong Q\oplus A$. \end{prop}

\begin{proof} By hypothesis, there is at least one infinite
idempotent $e\in R$, whence $eR$ is a directly infinite, finitely
generated projective right $R$-module. If $Q$ is any nonzero
finitely generated projective right $R$-module, Lemma
\ref{pinflem2} implies that
$Q$ is isomorphic to a direct summand of $eR$, and so $Q$ is
isomorphic to some nonzero right ideal $I\subseteq R$. Now since
$R$ is purely infinite, there is an infinite idempotent $f\in I$,
whence $fR$ is a directly infinite direct summand of $I$. It
follows that $I$ is directly infinite, whence $Q$ is directly
infinite.

The final conclusion of the proposition now follows immediately from
Lemma \ref{pinflem2}.
\end{proof}

\begin{theor}\label{pinfthm} Let $R$ be a  simple ring. Then $R$ is
purely infinite if and only if
\begin{enumerate}
\item[{\rm (a)}] $R$ is not a division ring, and
\item[{\rm (b)}] For every nonzero element $a\in R$, there exist
elements $x,y\in R$ such that $xay=1$.
\end{enumerate}
\end{theor}

\begin{proof} $(\Longrightarrow)$: Assume that $R$ is purely
infinite. Then $R$ contains an infinite idempotent, and so $R$ cannot
be a division ring.

Now consider a nonzero element $a\in R$. By assumption, there exists
an infinite idempotent $e\in aR$, and then Proposition
\ref{pinfprop} implies that $eR\cong R\oplus A$ for some $A$. Now
by Lemma \ref{pinflem1}, there exist orthogonal idempotents $f,g\in
R$ such that $e=f+g$ and
$fR\cong R$. Then $f\sim 1$, and so there are elements $\alpha\in
fR$ and
$\beta\in Rf$ such that $\alpha\beta=f$ and $\beta\alpha=1$. Since
$f\in eR \subseteq aR$, we also have $f=ar$ for some $r\in R$.
Therefore
$$1= \beta\alpha\beta\alpha= \beta f\alpha= \beta a(r\alpha),$$
and
(b) is established.

$(\Longleftarrow)$: Now assume conditions (a) and (b), and consider
a nonzero right ideal $I\subseteq R$. Since $R$ is not a division
ring,
$I$ must contain a proper nonzero right ideal, say $J$. Choose a
nonzero element $a\in J$, and apply condition (b): there exist
$x,y\in R$ such that $xay=1$. Now $e:= ayx$ is an idempotent lying
in
$aR$, whence $e\in I$ and $e\ne 1$. Since $(eay)(xe)=e$ and
$(xe)(eay)=1$, we also have $e\sim 1$, and so $1$ is infinite. But
then $e$, being equivalent to $1$, is infinite too, and we have
proved that $R$ is purely infinite.
\end{proof}

Notice that Theorem \ref{pinfthm} implies that the definition of a
purely infinite simple ring given above is left-right symmetric
(this can also be shown by a direct argument). Also, we point out
that in view of condition (b) of the theorem, our definition agrees
with the definition in current use among C*-algebraists. (For the
equivalence of this definition with Cuntz's original definition, see
\cite[Proposition 6.11.5]{Black}.)

\begin{corol}\label{pinfMorita} The class of  purely infinite
simple rings is closed under Morita equivalence.
\end{corol}

\begin{proof} It suffices to show that if $R$ is a  purely infinite
simple ring, then $M_2(R)$ and $eRe$ have the same properties,
where $e$ is an arbitrary nonzero idempotent in $R$. It is clear
that $M_2(R)$ and $eRe$ are simple.

Because $R$ contains an infinite idempotent, it cannot be artinian.
But $R$ is Morita equivalent to $eRe$, and thus $eRe$ cannot be a
division ring. For any nonzero element $a\in eRe$, Theorem
\ref{pinfthm} provides elements $x,y\in R$ such that $xay=1$, and
then
$exe$ and
$eye$ are elements of $eRe$ such that $(exe)a(eye)=e$. Thus, by the
theorem, $eRe$ is purely infinite.

By Proposition \ref{pinfprop}, there is an idempotent $f\in R$ such
that $fR
\cong 2R$, whence $fRf \cong M_2(R)$. Since $fRf$ is purely
infinite by the previous paragraph, we therefore conclude that
$M_2(R)$ is purely infinite, as desired.
\end{proof}

\begin{remas}
\label{pinfrem1} $\mbox{ }$ \vspace{.1truecm}

\begin{enumerate}
\item[(a)] Zhang (\cite[Theorem 1]{Zhang}), and also Brown and Pedersen
(\cite[Proposition 3.9]{BP}) have proved that every purely infinite
simple C*-algebra has real rank zero, which by \cite[Theorem
7.2]{agop} is equivalent to the property of being an exchange ring
\cite{Warfield}. Thus, the following question imposes itself:
\medskip

{\em Is every purely infinite simple ring an exchange ring?}
\medskip

\item[(b)] Let $R$ be any directly infinite, simple exchange ring
which is separative (see \cite{agop}). Then $R$ is necessarily
purely infinite, as follows. It follows from simplicity and the
exchange criterion of Nicholson and Goodearl that any nonzero right
ideal of $R$ contains a nonzero idempotent, say $e$. Then $R_R
\lesssim n(eR)$ for some $n$, whence $R$ is isomorphic to a
corner of $M_n(eRe)$, and so $M_n(eRe)$ is directly
infinite. By \cite[Proposition 2.3]{agop}, $eRe$ must be directly
infinite, that is, the idempotent $e$ is infinite.
\item[(c)] If $A$ is a simple exchange ring, then the simple ring $B=
\varinjlim M_{2^n}(A)$ (with block diagonal transition maps) will
either have stable rank $1$ or be purely infinite. For, notice that,
if there is some $n$ such that $M_{2^n}(A)$ is directly infinite,
then $B$ is purely infinite by Example \ref{pinfexamples}(e).
Otherwise,
$A$ is stably finite, and so it has power-cancellation by
\cite[Proposition 2.1.8]{RatC*}. This in turn implies that $B$ has
cancellation, and thus has stable rank $1$ by \cite[Theorem 9]{Yu}.
\item[(d)] If $A$ is a simple QB-ring (see \cite{ApedPer}), then
it either has stable rank $1$ or is purely infinite \cite{Perera}.
 \end{enumerate}
\end{remas}

\section{$K_0$ and $K_1$ of purely infinite simple rings}

Cuntz \cite{C2} computed the general patterns for $K_0$  and
(topological) $K_1$ of purely infinite simple C*-algebras. Here we
give parallel results for the $K_0$ and the (algebraic) $K_1$ groups
of purely infinite simple rings. We recall some basics of $K$-Theory
for the convenience of the reader.

Given a ring $R$, we define $\mathcal{V}(R)$ to be the set of
isomorphism classes (denoted $[A]$) of finitely generated projective
right $R$-modules, and we endow $\mathcal{V}(R)$ with the structure
of a commutative monoid by imposing the operation $$[A]+ [B] :=
[A\oplus B]$$ for any isomorphism classes $[A]$ and $[B]$. (The
notation $[A]$ will also be used for {\em stable} isomorphism
classes, that is, elements of $K_0(R)^+$, but it should be clear from
the context which is meant.) Equivalently
\cite[Chapter 3]{Black},
$\mathcal{V}(R)$ can be viewed as the set of equivalence classes of
idempotents in
$M_\infty(R)= \bigcup_{n=1}^\infty M_n(R)$ with the operation
$$[e]+[f] :=
\bigl[ \left( \smallmatrix e&0\\ 0&f
\endsmallmatrix \right) \bigr]$$
for idempotents $e,f\in
M_\infty(R)$. The group $K_0(R)$ is the universal group of $\mathcal{V}(R)$,
and the image of the canonical monoid homomorphism $\mathcal{V}(R)
\rightarrow K_0(R)$ is the {\em positive cone\/} of $K_0(R)$,
denoted $K_0(R)^+$.

Observe that $\mathcal{V}(R)$ is {\it conical\/}, that is, $x+y=0$ in $\mathcal{V}(R)$
only if $x=y=0$. Consequently, the subset $\mathcal{V}(R)^* := \mathcal{V}(R) \setminus
\{0\}$ is closed under addition. There is a natural pre-order $\le$
on $\mathcal{V}(R)$, where $x\le y$ if and only if there exists $z\in \mathcal{V}(R)$
such that $x+z=y$. In terms of isomorphism classes of finitely
generated projective modules $A$ and $B$, we have $[A] \le [B]$ if
and only if $A$ is isomorphic to a direct summand of $B$.  We say
that $\mathcal{V}(R)$ is {\it simple\/} provided that for each nonzero $x,y\in
\mathcal{V}(R)$, there exists $n\in\N$ such that $y\le nx$. Observe that if
$R$ is a  simple ring, then $\mathcal{V}(R)$ is a simple monoid (recall that
simplicity of $R$ implies that all nonzero projective $R$-modules
are generators in $\text{Mod-}R$). The converse is false (e.g.,
take $R=
\Z$).

\begin{lem}\label{pinflem3} Let $R$ be a  purely infinite simple
ring. For any $x,y\in \mathcal{V}(R)^*$, there exist $a,b\in \mathcal{V}(R)^*$ such that
$x=y+a$ and $y=x+b$.
\end{lem}

\begin{proof} This is just a restatement of the last part of
Proposition \ref{pinfprop}.
\end{proof}

\begin{prop} \label{vrgroup} If $R$ is a  purely infinite simple ring, then
$\mathcal{V}(R)^*$ is a group.

In particular, any nonzero finitely generated projective $R$-modules
which are stably isomorphic must be isomorphic.
\end{prop}

\begin{proof} According to \cite[Proposition 2.4]{agop}, to prove that
$\mathcal{V}(R)^*$ is a group it suffices to show that $\mathcal{V}(R)$
is conical and simple, and that for each $x\in
\mathcal{V}(R)^*$ there exists $b\in \mathcal{V}(R)^*$ such that $x+b=x$. However, we
already know that $\mathcal{V}(R)$ is conical and simple, and the remaining
condition follows from Lemma \ref{pinflem3} (take $y=x$). Therefore
$\mathcal{V}(R)^*$ is a group.

The final statement of the proposition now follows.
\end{proof}

An interesting consequence of Proposition \ref{vrgroup} is that if $R$
is a purely infinite simple ring and $[P]$ is the identity element of
$\mathcal{V}(R)^*$, then $P$ is nonzero and $P\oplus Q\cong Q$ for all
finitely generated projective right $R$-modules $Q$.

Notice that, aside from the fact that the stable rank of a
 purely infinite simple ring is $\infty$, Proposition \ref{vrgroup}
shows that the cancellative behavior of finitely generated projective
modules is almost the same as in the case of rings of stable
rank $1$.

If $M$ is an additive monoid, we write $\{0\} \sqcup M$ for the
monoid constructed by adjoining a new zero element to $M$. (To work
with this new monoid, one must choose some notation to distinguish
between the old and new zero elements.)

\begin{corol}\label{pinfVK0} If $R$ is a  purely infinite simple
ring, then $\mathcal{V}(R) \cong \{0\} \sqcup K_0(R)$, and $K_0(R)^+= K_0(R)$.
\end{corol}

\begin{proof} Let $\phi: \mathcal{V}(R) \rightarrow K_0(R)$ be the natural
monoid homomorphism; thus if $A$ is any finitely generated projective
right $R$-module, $\phi([A])= [A]$ is the stable isomorphism class
of $A$. By Proposition \ref{vrgroup}, $\mathcal{V}(R)^*$ is a group,
and so $\phi$ restricts to a group homomorphism $\phi^*:
\mathcal{V}(R)^* \rightarrow K_0(R)$. Thus the image of $\phi^*$ is a
subgroup of $K_0(R)$. In particular, this image contains $[0]$, and
so $\phi^*(\mathcal{V}(R)^*)= \phi(\mathcal{V}(R))= K_0(R)^+$. Since
$K_0(R)$ is generated by $K_0(R)^+$, we now see that $K_0(R)^+=
K_0(R)$.

If $x,y\in \mathcal{V}(R)^*$ and $\phi^*(x)= \phi^*(y)$, then $x+z=y+z$ for
some $z\in \mathcal{V}(R)$. Since either $z=0$ or $z$ lies in the group
$\mathcal{V}(R)^*$, it follows that $x=y$. This shows that $\phi^*$ is an
isomorphism of $\mathcal{V}(R)^*$ onto $K_0(R)^+= K_0(R)$, and therefore we
conclude that $\mathcal{V}(R) \cong \{0\} \sqcup K_0(R)$.
\end{proof}

For completeness, we present the following theorem which shows the
parallelism of patterns between algebraic and topological $K_1$ in the
cases of  purely infinite simple rings and C*-algebras. However, $K_1$
will not appear elsewhere in the paper.

\begin{theor} If $R$ is a  purely infinite simple ring then
$K_1(R)=U(R)^{\rm ab}$.
\end{theor}

\begin{proof} By \cite[Remark after 2.3]{memo2}, $R$ is a $GE$-ring,
so the natural map
$\kappa: U(R)\rightarrow K_1(R)$ is surjective. In order to show that
$\kappa$ is injective, we will proceed in two steps. The first one is
similar to Cuntz's argument in \cite{C2}.

{\it Step 1}. Let $v$ be a unit in the kernel of $\kappa$. Assume there
is a nonzero idempotent
$e\in R$ such that
$v=e+(1-e)v(1-e)$. Then $v\in U(R)'$.

{\it Proof of Step 1:} Take an idempotent $f<e$ such that $f\sim e$.
Set $r_1=1-e+f$, and note that $r_1\sim 1$. Now let $r_2, r_3, \dots
$ be orthogonal idempotents such that $r_2+ \cdots +r_i\le e-f$  and
$r_i\sim 1$ for all $i\geq 2$. Notice that $(r_1+ \cdots
+r_i)R(r_1+\cdots +r_i)\cong M_i(R)$; this isomorphism may be chosen
so that its restriction to $r_1Rr_1$ sends
$$1-e \longmapsto \mbox{diag}(1-e,0,\dots ,0) \qquad\quad \mbox{and}
\qquad\quad f \longmapsto \mbox{diag}(e,0,\dots ,0).$$
 Hence, the element
$v_i := (1-e)v(1-e)+f+r_2+\cdots +r_i$ corresponds to
$\text{diag}(v,1,\dots ,1)$ under this isomorphism. Since the image of
$v$ in $K_1(R)$ is
$0$, there is some $n\ge 1$ such that $\text{diag}(v,1,\dots,1)\in
GL_n(R)'$. Consequently, $v_n \in U \bigl( (r_1+ \cdots
+r_n)R(r_1+\cdots +r_n) \bigr)'$, and so $v= v_n+e-f -(r_2+ \cdots
+r_n)$ lies in $U(R)'$, as desired.

{\it Step 2}. For any $u\in U(R)$ there is a nonzero idempotent
$e\in R$ and a unit of the form $v=e+(1-e)v(1-e)$ such that
$u\equiv v \pmod{U(R)'}$.

{\it Proof of Step 2:} Any decomposition $1=e_1+\cdots+e_n$ with
$e_1\sim
\cdots \sim e_{n-1}$ and $e_n\lesssim e_1$ gives rise to an
isomorphism $R\cong S$, where
$S$ is the ring of $n\times n$ matrices $A= (a_{ij})$ over
$T=e_1Re_1$ of the following form: all the entries $a_{in}$  are in
$Tf$, and all the entries $a_{nj}$ are in $fT$, where $f$ is the
idempotent corresponding to $e_n$ under the subequivalence
$e_n\lesssim e_1$. An argument similar to that of Whitehead's Lemma
proves that, for $n\ge 4$, the set of elementary matrices relative
to the above matrix decomposition of $R$ is contained in $U(R)'$. Now
take
$n\ge 4$ and a decomposition of $1$ of the form indicated.
Consider the corresponding matricial representation over
$T=e_1Re_1$, which is also a purely infinite simple ring. By the
argument in
\cite[Theorem 2.2]{memo2}, the matrix corresponding to $u$ can be
transformed by elementary operations to a matrix of the form
$$\begin{pmatrix} 1& 0\\0& *\end{pmatrix} ,$$
where $*$ is a matrix of
size
$(n-1)\times (n-1)$. Since $n\ge 4$, the elementary matrices give
rise to elements in $U(R)'$, so that $u$ is congruent mod $U(R)'$
to a unit of the form
$e_1+(1-e_1)v(1-e_1)$, as desired.
\end{proof}

\section{Skew polynomial rings with freely independent indeterminates}
\label{skewpolysection}

In this section, we introduce a general construction that turns out
to be a keystone for our work. Let $R$ be a ring and $(\tau
_1,\delta _1),\dots ,(\tau_n,\delta _n)$ some right skew
derivations on $R$. Thus, the $\tau_i$ are ring endomorphisms of $R$,
and the $\delta_i$ are additive endomorphisms satisfying the rule
$\delta(rs)= \delta(r)\tau(s) +r\delta(s)$. Let $Y=\{y_1,\dots
,y_n\}$ be an $n$-element alphabet and denote by $Y^*$ the free monoid
on $Y$, with identity element denoted $1$. Label words in $Y^*$ in the
form $y_I=y_{i_1}y_{i_2}\cdots y_{i_t}$ for finite sequences
$I=(i_1,\dots ,i_t)$ of indices from $\{1,\dots ,n\}$, and similarly
for other quantities like $\tau_I$ or $\delta _I$.

We would like to build a ring whose right $R$-module structure is
free with basis $Y^*$, such that $ry_i=y_i\tau _i(r)+\delta _i(r)$
for all $i$ and all $r\in R$. To define a multiplication on this
$R$-module and verify the ring axioms is tedious; instead, we build
a ring of operators on this module, and afterwards carry over the
structure.

Let $V$ be the free right $R$-module with basis $Y^*$, and let
$E=\End_\Z (V)$. We let maps in $E$ act on the {\it right\/} of their
arguments, so that we can identify
$R$ with the subring of right multiplication operators in $E$.
Define
$z_1,\dots ,z_n\in E$ so that
$$(yr)z_i=(yy_i)\tau _i(r)+y\delta _i(r)$$
for all $i$, all $r\in
R$, and all $y\in Y^*$. Let $S$ be the subring of $E$ generated by
$R$ and $z_1,\dots ,z_n$. For all $i$, all $r,s\in R$, and all $y\in
Y^*$, we have
\begin{align} (ys)(rz_i) &= (ysr)z_i =(yy_i)\tau _i(sr)+y\delta
_i(sr) \notag\\
 &=(yy_i)\tau_i(s)\tau_i(r)+y\delta_i(s)\tau_i(r)+(ys)\delta_i(r)
=(ys)[z_i\tau_i(r)+\delta_i(r)]. \notag
\end{align}
Therefore
\begin{equation} rz_i=z_i\tau_i(r)+\delta_i(r) \tag{*}
\end{equation}
for all $i$ and all
$r\in R$. In particular, (*) implies that $S$ is generated as a right
$R$-module by the monomials $z_I$. All $\tau _i(1)=1$ and
$\delta_i(1)=0$, so
$(y)z_i=yy_i$ for all $y\in Y^*$. Consequently, we get by induction
that $(1)z_I=y_I$ for all $I$. Any $s\in S$ can be written as $\sum
_I z_Ir_I$ with almost all $r_I=0$, and $(1)s=\sum _I y_Ir_I$. Since
the $y_I$ form a basis for $V$ as a right $R$-module, the rule
$s\mapsto (1)s$ thus gives a right $R$-module isomorphism from $S$
onto $V$. We summarize our observations as follows:

\begin{prop}
\label{skew} Given a ring $R$ with right skew derivations $(\tau
_i,\delta _i)$ for $i=1,\dots ,n$, there is a ring $S$ containing
$R$ as a subring such that:

{\rm (a)} $S_R$ is free with basis $Y^*$, where $Y^*$ is the free
monoid on an $n$-element set $Y=\{y_1,\dots ,y_n \}$.

{\rm (b)} The ring and module multiplications $S\times R\rightarrow
S$ coincide.

{\rm (c)} The ring and monoid multiplications $Y^*\times
Y^*\rightarrow S$ coincide.

{\rm (d)} $ry_i=y_i\tau _i(r)+\delta _i(r)$ for all $i$ and all
$r\in R$.
\end{prop}

\begin{nota} We denote the ring $S$ in Proposition \ref{skew} by
$R\langle Y;\tau,\delta\rangle$, where $\tau$ and $\delta$ are
abbreviations for the $n$-tuples $(\tau _1,\dots ,\tau_n)$ and
$(\delta _1,\dots ,\delta_n)$. To see that $R\langle
Y;\tau,\delta\rangle$ is unique in a suitable sense, we establish
the following universal property:
\end{nota}

\begin{prop}
\label{univ} Let $R$ be a ring and $(\tau _1,\delta _1),\dots
,(\tau_n,\delta _n)$ right skew derivations on $R$. Suppose $\phi
:R\rightarrow T$ is a ring homomorphism and $t_1,\dots ,t_n\in T$
are elements such that $\phi (r)t_i=t_i\phi \tau_i(r)+\phi \delta_i(r)$
for all
$i$ and all
$r\in R$. Then there exists a unique ring homomorphism
$\ol{\phi}:R\langle Y;\tau,\delta\rangle\rightarrow T$ such that
$\ol{\phi}|_R=\phi$ and $\ol{\phi}(y_i)=t_i$ for all $i$.
\end{prop}

\begin{proof} It suffices to construct a ring with this universal
property and show it is isomorphic to $R\langle Y;\tau,\delta\rangle$.
More precisely, let $F=\Z\langle z_1,\dots ,z_n\rangle$ be the free ring
on $n$ letters, let $S_0=F*_{\Z}R$ be the ring coproduct of $F$ and $R$,
and let $S_1$ be the factor ring of $S_0$ by the ideal generated by
$(1*r)(z_i*1)-z_i*\tau _i(r)-1*\delta _i(r)$ for all $i$ and all $r\in
R$. Now let $\psi :R\rightarrow S_1$ be the map defined by $\psi
(r)=\ol{1*r}$ and set $w_i=\ol{z_i*1}\in S_1$.

Then for all $i$ and all $r\in R$ we have
\begin{equation}
\psi(r)w_i=w_i\psi\tau_i(r)+\psi\delta_i(r), \tag{$\dagger$}
\end{equation}
 and $(S_1,\psi,w_1,\dots ,w_n)$ is universal with respect to
$(\dagger)$. In particular, there is a unique ring homomorphism $\theta :S_1\rightarrow
R\langle Y;\tau,\delta\rangle$ such that $\theta \psi$ is the inclusion
map $R\rightarrow R\langle Y;\tau,\delta\rangle$ and $\theta (w_i)=y_i$
for all $i$. It is enough to show that $\theta$ is an isomorphism.

From $(\dagger)$, we see that $S_1=\sum _I w_I\psi (R)$. Note that
$\theta (\sum _I w_I\psi (r_I))=\sum _I y_Ir_I$ for all finite sums with
$r_I\in R$. Since $R\langle Y;\tau,\delta\rangle_R$ is free with basis
$\{y_I\}$, it follows that $\theta $ is an isomorphism.
\end{proof}

\begin{exem}\label{tyux} We give an example of the construction
considered above, and show that it coincides with a ring
constructed by Tyukavkin (\cite{tjufirst},
\cite{tju}).

Take
$R=k\langle\langle X\rangle \rangle$, where $X=\{x_1,\dots ,x_n\}$. Let
$\tau _1=\tau_2=\cdots =\tau_n=\tau$ be the unique $k$-algebra
homomorphism
$\tau :R\rightarrow R$ sending all $x_i$ to $0$.
Write elements of $R$ as infinite sums $\sum_{w\in X^*} \lambda_ww$,
where
$\lambda_w\in k$ and $X^*$ is the free monoid on $X$. Define
$k$-linear maps $\delta _i :R\rightarrow R$ by the rule
 $$\delta _i \bigl( \sum_{w\in
X^*} \lambda_ww \bigr)=\sum _{w\in X^*}\lambda_{wx_i}w.$$ We want to
check that the $\delta _i$ are right $\tau _i$-derivations. Take two
elements $r=\sum _w \lambda_ww$ and $r'=\sum _w \mu_ww$ in $R$, so that
$rr'=\sum _w \bigl( \sum _{uv=w} \lambda_u \mu_v \bigr)w$. Then we have
$$\delta_i(rr')= \sum _w \bigl( \sum _{uv=wx_i} \lambda_u \mu_v \bigr)w=
\sum _w  \bigl( \lambda_{wx_i}\mu_1+ \sum _{uv=w} \lambda_u \mu_{vx_i}
\bigr)w,$$ since $(wx_i)1$ is the only factorization $uv=wx_i$ where $x_i$
is not a right factor of $v$. Also,
\begin{align}
\delta_i(r)\tau (r')+r\delta_i(r') &= \bigl( \sum
_w \lambda_{wx_i}w \bigr)\mu_1 + \bigl( \sum _u \lambda_uu
\bigr) \bigl( \sum _v \mu_{vx_i}v \bigr)
\notag\\
 &=\sum_w \bigl( \lambda_{wx_i}\mu_1+\sum _{uv=w}
\lambda_u \mu_{vx_i} \bigr)w, \notag
\end{align}
 the same as above. This proves that $\delta _i$ is a right
$\tau _i$-derivation. It follows that there exists a $k$-algebra
$R\langle Y;\tau,\delta\rangle$. Note that $\tau _i(x_j)=0$ and
$\delta _i(x_j)=\delta _{ij}$ for all $i,j$. Hence, $x_iy_j=
\delta_{ij}$ for all $i,j$. It follows that the algebra $R\langle
Y;\tau,\delta\rangle$ in this case coincides with the algebra
constructed by Tyukavkin in \cite[page 404]{tju}.
\end{exem}

The examples of Leavitt also arise from the $R\langle
Y;\tau,\delta\rangle$
construction, as we observe in the next section.

\section{Leavitt's algebras}
\label{first}

This section is devoted to showing that some universal examples of
non-IBN algebras lie in the class of purely infinite simple rings.
Our interest in quoting them here is twofold: on one side, these
examples are the algebraic precursors of Cuntz algebras; on the other
side, Tyukavkin's example quoted in Example \ref{tyux} are, in some
sense, completions of Leavitt's examples.

For any field $k$, and for any two natural numbers $m,n$, Leavitt
(\cite{Leavitt}) introduced the $k$-algebras $V_{m,n}$, with a
universal isomorphism $i: nV_{m,n}\longrightarrow mV_{m,n}$, and
$U_{m,n}$, with a universal pair of morphisms $i:
nU_{m,n}\longrightarrow mU_{m,n}$ and $j: mU_{m,n}\longrightarrow
nU_{m,n}$ such that $ji=1_{nU_{m,n}}$. Later, Cohn, Skornyakov and
Bergman proved some fundamental properties of these algebras. The
monoid $\mathcal{V}(R)$ of these examples was computed by Bergman in
\cite{Bergman}. For the sake of completeness, we recall here
Bergman's result.

\begin{theor} \label{uniBerg} {\rm (\cite[Theorem 6.1] {Bergman})}
For a field $k$ and positive integers $m,n$, the rings $V_{m,n}$ and
$U_{m,n}$ are hereditary. Moreover, $\mathcal{V}(V_{m,n})=\langle
I\mid mI=nI\rangle$, and  $\mathcal{V}(U_{m,n})=\langle I, P\mid
mI=nI+P\rangle$.
\end{theor}

We will use the constructions in Section 3 to give some extra
information on the algebras $V_{1,n}$ and $U_{1,n}$. First of all,
notice that
$U_{1,n}= k\langle X\rangle \langle Y;\tau,\delta \rangle$, where
$X$, $Y$, $\tau$, $\delta$ are as in Example \ref{tyux}. This is clear
once we observe that the row $(y_1,\dots ,y_n)$ and the column
$(x_1,\dots ,x_n)^T$ give a universal pair of morphisms $i:
nR_n\longrightarrow R_n$ and $j: R_n\longrightarrow nR_n$ respectively
such that
$ji=1_{nR_n}$, where we set $R_n=k\langle X\rangle \langle
Y;\tau,\delta \rangle$. The algebra $V_{1,n}$ is thus the algebra
obtained by imposing the relation $ij=1_{R_n}$ in $R_n$, that is,
$V_{1,n}\cong R_n/I_n$, where $I_n$ is the ideal of $R_n$ generated
by the idempotent $e_n=1-\sum _{i=1}^n y_ix_i$. Following Tyukavkin
(\cite{tjufirst}, \cite{tju}) we will call the elements in $X^*$
{\em monomials\/} and the elements in $Y^*$ {\em words\/}. From our
basic relations
$x_iy_j=\delta _{ij}$ it follows that to each monomial $x_I$ there
corresponds a word $y_{I^*}$ such that $x_Iy_{I^*}=1$, where $I^*
=(i_r,\dots ,i_1)$ for $I=(i_1,\dots ,i_r)$. A {\it monoword} is any
element of the form $y_Ix_J$ for some indices $I,J$. Notice that the
monowords $y_Ix_J$ form a $k$-basis of $U_{1,n}$.

\begin{theor} \label{otro}
For every natural number $n\geq 2$, $V_{1,n}$ is a purely infinite
simple ring, and $K_0(V_{1,n}) \cong \Z/(n-1)\Z$.
\end{theor}

\begin{proof} By Theorem 4.1, $\mathcal{V}(V_{1,n})= \langle I\mid
nI=I \rangle$, from which it follows that $K_0(V_{1,n}) \cong
\Z/(n-1)\Z$.

We identify $k\la X\ra\la Y;\tau ,\delta \ra$ with
$U_{1,n}$. Let $I_n$ be the ideal of $U_{1,n}$ generated by $e_n := 1-
\sum_{i=1}^n y_ix_i$. As observed before, we have $V_{1,n}\cong
U_{1,n}/I_n$. Every element $\alpha\in U_{1,n}$ can be written
uniquely as $\alpha= \sum_{I,J} \lambda_{I,J}y_Ix_J$ for scalars
$\lambda_{I,J}$ which are almost all zero. The {\em support\/} of
$\alpha$ is the set of pairs $(I,J)$ such that $\lambda_{I,J} \neq 0$;
denote by $s(\alpha)$ the cardinality of the support of $\alpha$.
Since $x_ie_n= e_ny_j=0$ for all $i,j$, the ideal $I_n$ is spanned
(over $k$) by products of the form $y_Ie_nx_J$. It follows that the
support of any nonzero element of $I_n$ must contain a pair $(I,J)$
with both $I$ and $J$ nonempty. Consequently, $I_n$ contains no
nonzero elements of either $k\la X\ra$ or $k\la Y\ra$.

We
prove that for $\alpha \in U_{1,n}\setminus I_n$, there exist
$\beta,\gamma
\in U_{1,n}$ such that $\beta\alpha\gamma =1$. (In fact, $\beta$ can
be taken as a scalar times a monomial, and $\gamma$ as a word.)
If $s(\alpha)=1$, then $\alpha= \lambda_{I,J}y_Ix_J$ with
$\lambda_{I,J} \neq 0$, and so
$$(\lambda_{I,J}^{-1} x_{I^*})
\alpha y_{J^*} =1.$$
Now assume that $s(\alpha) =d>1$. There must be
an index $i$ such that $\alpha y_i \notin I_n$, since otherwise we
would have $\alpha= \alpha \bigl( e_n+ \sum_{i=1}^n y_ix_i \bigr) \in
I_n$. Note that either $\alpha \in k\langle Y\rangle$ or the total
degree in $X$ of $\alpha y_i$ is less than that of $\alpha$. Hence,
by induction on this degree, there is a word $y_K$ such that $\alpha
y_K$ is a nonzero polynomial in $k\langle Y\rangle$. Clearly $\alpha
y_K \notin I_n$, and so there is a monomial $x_I$ such that $\alpha'=
x_I\alpha y_K$ is a polynomial in $k\langle X\rangle$ with nonzero
constant term. Clearly, we can choose a word $y_J$ such that
$\alpha'y_J \neq 0$ and $s(\alpha'y_J) <d$. Since $x_{J^*}\alpha'
y_J$ is a polynomial in $k\langle X\rangle$ with nonzero constant
term, we see that $\alpha' y_J \notin I_n$. By induction, there
exist $\beta,\gamma \in U_{1,n}$ such that $\beta x_I\alpha y_K y_J
\gamma= \beta\alpha' y_J\gamma =1$. This establishes the claim.

By Theorem \ref{pinfthm}, since $V_{1,n}$ is clearly not a division
ring, it is simple and purely infinite.
\end{proof}

Notice that we may identify
each $U_{1,n}$ with $k\langle X_n\rangle \langle Y_n; \tau, \delta
\rangle$ where $X_n= \{x_1,\dots,x_n\}$ and $Y_n= \{y_1,\dots,y_n\}$
are subsets of fixed infinite sets $\{x_1,x_2,\dots \}$ and $\{y_1,
y_2,\dots \}$. In particular,
$U_{1,n} \subset U_{1,n+1}$ for all $n$, and we set $U_\infty$ equal to
the union of the $U_{1,n}$.

\begin{theor} \label{inflim1} The ring $U_{\infty }$ is simple and
purely infinite. Moreover, $K_0(U_{\infty })\cong \Z$.
\end{theor}

\begin{proof} Clearly, $U_{\infty}$ is not a division ring, because
it contains at least one infinite idempotent.

Let $\alpha$ be a nonzero element of $U_\infty$. Then $\alpha \in
U_{1,n}$ for some $n$, and $\alpha
=\sum_{I,J}{\lambda _{I,J}y_Ix_J}$ where the $\lambda_{I,J}
\in k$, the $y_I$ are monomials in $y_1,\dots,y_n$, and the $x_J$ are
words in $x_1,\dots,x_n$. We choose an index $I'$ of minimal length
with
$\lambda _{I',J}\ne 0$ for some $J$, and then, we choose an index
$J'$ of maximal length among the indices $J$ such that $(I',J)$ is in
the support of $\alpha$. Then,
$$({\lambda _{I',J'}}^{-1}x_{{I'}^*})
\alpha ( y_{{J'}^*})=1+\sum_{I,J}{\lambda '
_{I,J}y_{I}x_{J}}\, ,$$
and $\lambda'_{\emptyset,J} =0$ for all $J\ne \emptyset$. Now, since
$x_{n+1}, y_{n+1}
\in U_{\infty }$ and $x_{n+1}y_I=0$ for all nontrivial words $y_I$ in
$U_{1,n}$, we have $(x_{n+1}{\lambda _{I',J'}}^{-1}x_{{I'}^*}) \alpha
( y_{{J'}^*}y_{n+1})=1$. Thus $U_{\infty }$ is simple and purely
infinite because of Theorem \ref{pinfthm}.

By Theorem \ref{uniBerg}, we have $K_0(U_{1,n})\cong \Z$ and a
generator of $K_0(U_{1,n})$ is provided by the class
$[U_{1,n}]$. It follows that the induced
homomorphisms $K_0(U_{1,n}) \rightarrow
K_0(U_{1,n+1})$ are isomorphisms. Since the functor $K_0(-)$
preserves direct limits, we get $K_0(U_{\infty })\cong \varinjlim
K_0(U_{1,n}) \cong \Z.$
\end{proof}

\section{Algebras $R\langle Y; \tau,\delta \rangle$ with other
coefficients}
\label{SandU}

Now we proceed to exploit the skew polynomial construction of Section
\ref{skewpolysection}, in order to get new examples of purely infinite
simple rings. These will be, in some sense, ``intermediate'' algebras
between Leavitt's and Tyukavkin's examples, and, as we will see, they
all  have the same K-theoretical behavior. We begin
as in Example
\ref{tyux}, except that we now index our variables starting at 0
rather than at 1. For a given natural number
$n$, consider the algebra
$k\langle\langle X_n\rangle \rangle$, where $X_n=\{x_0,x_1,\dots
,x_n\}$, equipped with the skew derivations $(\tau_i,\delta_i)$
defined in Example \ref{tyux}. Throughout this section, $R_n$ will
denote a local subalgebra of
$k\la\la X_n\ra\ra$ containing $k\la X_n \ra$ such that $R_n$ is
invariant under the $\delta_i$. Of course, one possibility is $R_n=
k\langle\langle X_n\rangle \rangle$; another choice of $R_n$ will be
important in Section \ref{isoQT}. Since the augmentation ideal of
$k\langle\langle X_n\rangle \rangle$ has codimension $1$, its
intersection with $R_n$ must be the maximal ideal of $R_n$. Thus, all
power series in $R_n$ with nonzero constant term are invertible in
$R_n$. Note that by \cite[Proposition 2.9.18]{free}, $R_n$ is a
semifir.

By assumption, the $(\tau_i,\delta_i)$ restrict to skew derivations
on $R_n$. Thus, we may
construct the skew polynomial algebra $S_n=R_n\la
Y_n;\tau,\delta \ra$, where $Y_n= \{y_0,y_1,\dots,y_n\}$, and our first
goal is to determine the structure of this algebra. Except in
Proposition \ref{ibn2}, where we consider the sequence of algebras
$S_1,S_2,\dots$, we shall throughout the remainder of the section keep
$n$ fixed and write $R$, $S$, $X$, $Y$ for $R_n$, $S_n$, $X_n$, $Y_n$.

\begin{lem}
\label{rwr'} {\rm (cf.~\cite{tju})} For any nonzero element $r\in R$
there are $w\in Y^*$ and $r'\in R$ such that $rwr'=1$. In fact,
given any nonzero elements $r_1,\dots,r_m \in R$, there exists
$w\in Y^*$ such that $r_jw\in R$ for all $j$ and some $r_iw$ is
invertible in
$R$.
\end{lem}

\begin{proof} Let $r_1,\dots,r_m$ be nonzero elements of $R$. For
$j=1,\dots,m$, let
$l_j$ denote the minimum length of monomials occurring in $r_j$. We
may assume that $l_1\leq \dots\leq l_m$. Pick a monomial $x_I$ of
length $l_1$ occurring in $r_1$, and set $w= y_{I^*}$. Observe that
for any monomial $x_J$ of length at least $l_1$, either $x_Jw=0$ or
$x_Jw$ is a monomial. Hence, $r_jw\in R$ for all $j$, and $r_1w$
has nonzero constant term. Consequently, $r_1w$ is invertible in
$R$, by our choice of
$R$.
\end{proof}

\begin{prop}
\label{soc} {\rm (a)} The algebra $S$ is a prime ring, and its socle is
the ideal $SeS$ generated by the minimal idempotent $e := 1-\sum
_{i=0}^n y_ix_i$. Further, $\soc (S)$ is a regular ideal of $S$.

{\rm (b)} Suppose that $\phi:R\rightarrow V$ is a $k$-algebra
homomorphism and that there exist elements $t_0,\dots,t_n \in V$ such
that $\phi(x_i)t_j= \delta_{ij}$ for all $i,j$. Then $\phi$ extends
uniquely to a $k$-algebra homomorphism $\ol{\phi} :S\rightarrow V$
such that $\ol{\phi}(y_j)= t_j$ for all $j$.
\end{prop}

\begin{proof} (a) Note that $ey_j= x_je= 0$ for all $j$. Consequently,
$eS= eR$ and $eSe=ke$.

 Given a nonzero
element $\alpha $ in $S$, either $e\alpha =\alpha$ or there is some
$i$ such that $x_i\alpha \ne 0$. It follows by induction on the
maximum length of words occurring in $\alpha$ that there is a
monomial $m\in X^*$ such that $em\alpha \ne 0$. By the observation
above, $em\alpha =er$ for some $r\in R$, and $r$ is right invertible
in $S$ by Lemma \ref{rwr'}. In particular, there is some $\gamma\in
S$ such that $em\alpha\gamma =e$. It follows at once that $S$ is a
prime ring and that $SeS$ is the unique minimal nonzero ideal of $S$.
Then, because $eSe= ke\cong k$ we conclude that $eS$ is a minimal
right ideal of $S$. Therefore $SeS= \soc (S)$.

 Now
$SeS$ is a simple ring having a minimal one-sided ideal, and
Litoff's Theorem (see \cite{FU}) says that it is locally a matrix
ring over a division ring (actually the division ring is $k$ in our
case). In particular
$\soc (S) =SeS$ is a regular ring.

(b) In view of Proposition 3.3, it suffices to show that
\begin{equation}
\phi(r)t_j= t_j\phi\tau(r)+ \phi\delta_j(r)
\end{equation}
for any $r\in R$ and all $j=0,\dots,n$. We can write $r= \alpha+
r_0x_0+ \dots+ r_nx_n$ for some $\alpha\in k$ and $r_i\in
k\langle\langle X\rangle\rangle$. Then $\alpha= \tau(r)$ and $r_i=
\delta_i(r)$; in particular, all the $r_i \in R$. Thus
$$\phi(r)t_j= \bigl( \alpha+ \sum_{i=0}^n \phi(r_i)\phi(x_i) \bigr)
t_j= \alpha t_j+ \phi(r_j)= t_j\phi\tau(r)+ \phi\delta_j(r),$$
proving (1).
\end{proof}

The next lemma may be known,
but we have not been able to locate any reference. Write $M^*= \Hom
_R(M,R)$ for $R$-modules $M$.

\begin{lem}
\label{ulreg} Let $R$ be a left semihereditary ring, and let
$T=R_{\Sigma}$ be a universal localization of $R$. Assume that for
all finitely presented right $R$-modules $M$ such that $M^*=0$, we
have $M\otimes _RT=0$. Then $T$ is a regular ring and
every finitely generated projective  right
$T$-module is induced from a finitely generated projective right
$R$-module.
\end{lem}

\begin{proof} We first claim that any finitely presented right
$R$-module $M$ has a decomposition $M= M_1\oplus M_2$ with $M_1$
projective and $M_2^*=0$. In case $R$ is semihereditary on both
sides, this follows from \cite[Theorem 1.2(3)]{Lueck} (the
hypothesis of an involution is not needed in the proof). Let
\begin{equation} \begin{CD}
R^m @>f>> R^n @>g>> M @>>>0 \end{CD} \notag
\end{equation}
be a resolution of $M$. Taking duals, we obtain an exact sequence
\begin{equation} \begin{CD}
0 @>>> M^* @>g^*>> {}^nR @>f^*>> {}^mR. \notag
\end{CD}\end{equation}
Now $f^*({}^nR)$ is a finitely generated submodule of $^mR$, so it is
projective because $R$ is left semihereditary. Consequently,
$g^*(M^*)= ({}^nR)e$ for some idempotent matrix $e\in M_n(R)$. It
follows that
$$(1-e)(R^n)= g^{-1} \bigl( \{x\in M\mid h(x)=0 \
{\mbox{\rm for all}}\ h\in M^*\} \bigr).$$
In particular, $\ker(g) \subseteq (1-e)(R^n)$. Hence, $M= M_1\oplus
M_2$ where $M_2= g \bigl( (1-e)(R^n) \big)$ satisfies $M_2^*=0$ and
$M_1 \cong e(R^n)$ is projective.

Now let $N$ be a finitely presented right $T$-module. By
\cite[Corollary 4.5]{scho}, there exists a finitely presented right
$R$-module $M$ such that $M\otimes _RT\cong N$. By the above, we can
write $M=M_1\oplus M_2$ with $M_1$
projective and $M_2^*=0$. By hypothesis, we
have
$M_2\otimes _RT=0$, and so $N\cong M_1\otimes_RT$ is an induced
finitely generated projective right $T$-module.
\end{proof}

\begin{theor}
\label{main} Let $R$, $S$, and $e$ be as above, and set
$I=SeS$. Let
$f:R\rightarrow R^{n+1}$ be the homomorphism given by left
multiplication by the column $(x_0,\dots, x_n)^T$. Then the universal
localization
$R_f$ is a  purely infinite simple regular ring, and every finitely
generated projective $R_f$-module is free. Moreover, $R_f\cong S/I$,
and $S$ is regular.
\end{theor}

\begin{proof} We first prove that $R_f\cong S/I$. For that, it is
enough to show that the natural map $R\rightarrow S/I$ satisfies the
universal property of the universal localization $R\rightarrow R_f$.
If $\phi :R\rightarrow V$ is a $k$-algebra homomorphism such that
$f\otimes 1_V$ is invertible, then there are elements $t_0,\dots
,t_n$ in $V$ such that  $\phi (x_i)t_j=\delta_{ij}$ and $\sum
_{i=0}^nt_i\phi (x_i)=1$. By Proposition \ref{soc}(b), there exists a
unique $k$-algebra map $\ol{\phi}:S\rightarrow V$ such that
$\ol{\phi}|_R=\phi$ and $\ol{\phi}(y_i)=t_i$ for all $i$. Clearly
this map factors through $S/I$ and so we get the desired map from
$S/I$ to $V$. This shows that we may identify $S/I$ with $R_f$. It
follows that the localization map $R\rightarrow R_f$ is injective
(recall from Lemma \ref{rwr'} that nonzero elements of $R$ are right
invertible in $S$).

Write $T=S/I=R_f$, and identify $R$ with its image in $T$. Let us
check directly that $T$ is purely infinite simple. (This also follows
from the fact that $T$ is regular with every finitely generated
projective $T$-module being free, which we will prove later.) Let
$\alpha$ be an element of $S$ which is not in $I$. There must exist
$i$ such that $x_i\alpha \notin I$, since otherwise $\alpha =e\alpha
+\sum_{i=0}^n y_ix_i\alpha \in I$. Note that either $\alpha \in R$ or
the degree in
$Y$ of
$x_i\alpha $ is smaller than that of $\alpha$. We conclude that there
is a monomial $m\in X^*$ such that $m\alpha$ is a nonzero element of
$R$. By Lemma \ref{rwr'}, we get an element $g\in S$ such that
$m\alpha g=1$. This shows that $T$ is a purely infinite simple ring.

Let $M$ be a finitely presented right $R$-module such that $M^*=0$.
We want to prove that $M\otimes _RT=0$, in order to apply Lemma
\ref{ulreg}. Take a presentation
$$R^s\rightarrow R^t\rightarrow
M\rightarrow 0.$$
Since the functor $(-)\otimes _RT$ is right exact,
we get an exact sequence
$$T^s\rightarrow T^t\rightarrow M\otimes
_RT\rightarrow 0.$$
Write $z_i$ for the image of $y_i$ in $T$. The
map $R^s\rightarrow R^t$ above is given by left multiplication by a
$t\times s$ matrix $A$ with coefficients in $R$, and we have to see
that $AT^s=T^t$, that is, the columns of $A$ generate $T^t$ as a right
$T$-module. We proceed  by induction on $t$. If $t=1$, then there
exists a nonzero entry $p$ in $A$, because $M^*=0$. By Lemma
\ref{rwr'}, the element $p$ is right invertible in $S$, so in $T$,
and thus the entries of $A$ generate $T$ as a right $T$-module. Now
assume that $t>1$. If some entry of $A$ is invertible in $R$, then by
a standard process we can find invertible matrices $P$ and $Q$ over $R$
of appropriate sizes such that $PAQ= \left(
\begin{array}{cc} 1 & 0 \\ 0 & A'\\
\end{array}\right)$, where $A'$ is a matrix of size $(t-1)\times
(s-1)$. In this case, $M$ can be generated by $t-1$ elements, and
induction applies.

In the general case, apply Lemma \ref{rwr'} to the entries $a_{ij}$
of $A$, to obtain a word $w\in Y^*$ such that all $a_{ij}w \in R$
with at least one $a_{ij}w$ invertible. Consequently, there is a
vector $v\in T^s$ of the form $(0,\dots,0,\ol{w}, 0,\dots,0)^T$
such that the column $q= Av\in AT^s \cap R^s$ has
some entry which is invertible in $R$. Now consider the matrix
$C=\left(\begin{array}{cc} q & A \\
\end{array}\right)$ with coefficients in $R$, of size $t\times
(s+1)$. The finitely presented module $M':=R^t/CR^{s+1}$ is a factor
module of $M$, and consequently $(M')^*=0$. Moreover, the matrix $C$
has an invertible entry, and so we get as before by induction that
$M'\otimes _RT=0$, or equivalently that $T^t=CT^{s+1}$. But $q$
belongs to the $T$-submodule of $T^t$ generated by the columns of
$A$, and therefore the columns of $A$ generate $T^t$ as an
$T$-module. We conclude that $M\otimes _RT=0$, as desired.

By \cite[Proposition 2.9.18]{free}, $R$ is a semifir. It follows
from Lemma \ref{ulreg} that $T$ is a regular ring and
that every finitely generated projective right $T$-module is free.

Finally, $S/I$ is a regular ring, and by Proposition \ref{soc}, $I$
is a regular ideal of $S$. Thus  it follows from \cite[Lemma
1.3]{vnrr} that $S$ is regular.
\end{proof}

The regularity of $S$ was proved by Tyukavkin
(\cite{tjufirst}, \cite{tju}) in the case where $R= k\langle\langle
X\rangle\rangle$. In order to compute
$K_0$ groups of both $S$ and $R_f$, we will need the following
technical lemma.

As in the proof of Theorem \ref{main}, we set $T= S/I= R_f$, and we
write
$T_n$ in case $n$ requires mention.

\begin{lem} \label{ibn1} If there is a right $R$-module map
$p:R\rightarrow R^s$ which becomes invertible over
$T$, then $n$ divides $s-1$.
\end{lem}

\begin{proof} Write $p=(p_1,\dots ,p_s)^T$, where each
$p_i\in R$. We will construct, by induction on $i$, words $w_i$ in
$Y^*$ and invertible elements
$g_i$ in $R$ such that the
following statements hold:
\begin{enumerate}
\item[$(P_i)$]  There exists an invertible map  $p^{(i)}:T\rightarrow
T^s$ satisfying the following properties:
\subitem(1) $p^{(i)}_{i+1},\dots ,p^{(i)}_{s}\in R$.
\subitem(2) The inverse of $p^{(i)}$ is the row $(w_1g_1,\dots
,w_ig_i,\alpha _{i+1},\dots ,\alpha _s )$ for some elements $\alpha
_{i+1},\dots ,\alpha_s\in T$.
\end{enumerate}

The statement is obvious when $i=0$. Assume that $0\le i<s$, and
that
$(P_i)$ holds. We will prove $(P_{i+1})$. Note that $p^{(i)}_jw_m=0$
for $j\ne m$. Without loss of generality, we can assume that the
order of the series $p^{(i)}_{i+1}$ is less than or equal to the
order of
$p^{(i)}_{i+t}$ for all $t\ge 2$. Choose a word $w_{i+1}\in Y^*$
with length equal to the order of $p^{(i)}_{i+1}$ and such that
$p^{(i)}_{i+1}w_{i+1}$ is invertible in $R$. Let $g_{i+1}\in R$ be
the inverse of $p^{(i)}_{i+1}w_{i+1}$ and note that
$$1=p^{(i)}_{i+1}w_{i+1}g_{i+1}=p^{(i)}_{i+1}\alpha _{i+1}.$$ It
follows that $u=\alpha
_{i+1}p^{(i)}_{i+1}+(1-w_{i+1}g_{i+1}p^{(i)}_{i+1})$ is invertible
in
$T$ with inverse $u^{-1}=w_{i+1}g_{i+1}p^{(i)}_{i+1}+(1-\alpha
_{i+1}p^{(i)}_{i+1})$. Therefore $p^{(i+1)}:=p^{(i)}u$ is invertible
with inverse
$$u^{-1}(w_1g_1,\dots ,w_ig_i,\alpha _{i+1},\dots ,\alpha _s).$$

Note that, for $t>1$, we have
$$p^{(i)}_{i+t}u=p^{(i)}_{i+t}(1-w_{i+1}g_{i+1}p^{(i)}_{i+1}).$$
Since the order of $p^{(i)}_{i+t}$ is greater than or equal to the
length of $w_{i+1}$, we conclude that $p^{(i+1)}_{i+t}\in R$, and
condition (1) of $(P_{i+1})$ holds. On the other hand, for $m\le i$
we have $p^{(i)}_{i+1}w_m=0$ and so
$$u^{-1}w_mg_m =w_mg_m.$$
We also
have $$u^{-1}\alpha _{i+1}= w_{i+1}g_{i+1}p^{(i)}_{i+1}\alpha
_{i+1}+(1-\alpha _{i+1}p^{(i)}_{i+1})\alpha _{i+1}=w_{i+1}g_{i+1}, $$
and so condition (2) of $(P_{i+1})$ is also satisfied. Therefore the
induction works.

Take $q_i=g_ip^{(s)}_i\in T$ for $i=1,\dots,s$. Then
\begin{equation}
\sum_{i=1}^s w_iq_i=1, \tag{1}
\end{equation}
$q_iw_i\ne 0$ for all $i$, and
$q_iw_j=0$ for $i\ne j$. We claim that these conditions imply
$s\equiv 1 \pmod{n}$. We proceed by induction on the maximum of
the lengths of the $w_i$. This maximum is $0$ if and only if $s=1$
(and then $q_1=1$). So assume that $s>1$. In this case all $w_i$ are
different from $1$. Fix $\ell \in\{0,\dots ,n\}$. Left multiplying (1)
by $x_\ell$ and right multiplying it by $y_\ell$, and letting
$A_\ell=\{i: x_\ell w_i\ne 0\}$, we have
$$\sum _{i\in A_\ell}(x_\ell w_i)(q_iy_\ell) =1.$$
Note that $\{1,\dots ,s\}$ is the disjoint union of the family
$\{A_\ell\mid
\ell=0,\dots ,n \}$. Observe also that for $i,j\in A_\ell$ we have
$(q_iy_\ell)(x_\ell w_j)=q_iw_j$. So this term is $0$ if $i\ne j$ and
nonzero if $i=j$. By induction, $|A_\ell|\equiv 1 \pmod{n}$. Therefore
$$s=\sum _{\ell=0}^n |A_\ell|\equiv n+1\equiv 1 \qquad (\mbox{\rm
mod\ } n),$$ as desired.
\end{proof}

\begin{theor}
\label{K0} Let $R$, $T$, and $f:R\rightarrow R^{n+1}$ be as above. Then
$K_0(T)$ is a cyclic group of order $n$, generated by $[T]$.
\end{theor}

\begin{proof}
Consider the homomorphism $\iota ^*:K_0(R)\rightarrow K_0(T)$, where
$\iota :R\rightarrow T$ denotes the localization map. Since our ring
$R$ is a semifir, we have $K_0(R)$ infinite cyclic, generated by
$[R]$. Also, by Lemma \ref{ulreg} and the proof of Theorem
\ref{main}, the map $\iota^*$ is a group epimorphism. Thus, it is
enough to show that the kernel of $\iota^*$ is generated by $n[R]$ to
get the desired result.

Suppose that $[R^r]-[R^s]$ lies in $\ker(\iota^*)$, for some
nonnegative integers $r$, $s$, that is, $[T^r]- [T^s] =0$ in
$K_0(T)$. Then $T^{r+t} \cong T^{s+t}$ for some $t\geq 0$. Since
$[R^r]-[R^s]= [R^{r+t}]- [R^{s+t}]$, we may thus assume that $T^r
\cong T^s$. By \cite[Corollary 4.4]{scho}, there is an isomorphism
$T^{r+u} \rightarrow T^{s+u}$, for some $u\geq 0$, which is induced
from a map $R^{r+u} \rightarrow R^{s+u}$. We can again reduce to the
case that $u=0$. Thus, it suffices to prove that if there exists a map
$g: R^r \rightarrow R^s$ that becomes invertible over
$T$, then $n$ divides $r-s$. After taking the direct sum of $g$
with an identity map
$R^m \rightarrow R^m$ for a suitable $m$, we may assume that $r=
\ell n+1$ for some positive integer $\ell$. Since the map
$f:R\rightarrow R^{n+1}$ becomes an isomorphism over $T$, there exists
a homomorphism $h: R\rightarrow R^r$ which becomes invertible over $T$,
and the same holds for $gh: R\rightarrow R^s$. By Lemma \ref{ibn1},
$n$ divides $s-1$, and hence also $r-s$, as
desired.
\end{proof}

\begin{corol}\label{KzeroS}
Let $R$ and $S$ be as above. For $\ell,m\in\N$, we have $S^\ell \cong
S^m$ if and only if $\ell=m$. Moreover,
$K_0(S)$ is an infinite cyclic group, with generator $[S]$.
\end{corol}

\begin{proof} Suppose that $S^\ell \cong
S^m$ for some $\ell\geq m>0$.
Consider $p\in \N$ such that $p\geq n$ and $p> \ell-m$. If we take
$R_p= k\la\la X_p\ra\ra$, then $S$ embeds in $S_p$, whence $S_p^\ell
\cong S_p^m$ and so $T_p^\ell \cong T_p^m$. Thus, by Theorem \ref{K0}
(applied to
$T_p$), we get that
$p$ divides $\ell-m$, which implies that $\ell=m$, as desired.

Now, by the above argument, we have that the class $[S]\in
K_0(S)$ generates an infinite cyclic subgroup. By Proposition
\ref{soc}, $I$ is the socle of $S$, and $I$ is generated by the
minimal right ideal $eS$. Hence, $K_0(I)$ is infinite cyclic, with
generator $[eS]$. Since $y_0x_0, \dots, y_nx_n$ are orthogonal
idempotents equivalent to $1$, we have $eS \oplus S^{n+1} \cong S$,
and so $[eS]= -n[S]$ in $K_0(S)$. In particular, it follows that the
natural map
$K_0(I)\rightarrow K_0(S)$ is injective. Using that, we get the
following commutative diagram, whose bottom row is exact because $S$
is regular.
\begin{equation}\begin{CD}
0 @>>> \Z @>{-n}>> \Z @>>> \Z/n\Z @>>> 0\\
@. @VV{1\mapsto[eS]}V @VV{1\mapsto[S]}V @VV{[1]_n\mapsto[T]}V @.\\
0 @>>> K_0(I) @>>> K_0(S) @>>> K_0(T) @>>> 0
\end{CD} \notag \end{equation}
Moreover, the map $\Z/n\Z \rightarrow K_0(T)$ is an
isomorphism by Theorem \ref{K0}. Thus, by the Five Lemma we conclude
that $K_0(S)$ is infinite cyclic with generator $[S]$, as desired.
\end{proof}

Now consider a sequence of algebras $S_1,S_2,\dots$ of the type
discussed in this section. Assume the sequence has been chosen so that
$R_1\subset R_2\subset \dots$; we could, of course, take $R_n= k\la\la
X_n\ra\ra$ for all $n$, but different choices will be needed in
Section 8. We have natural inclusions
$S_n
\subset S_{n+1}$ for all
$n$, and we set $S_{\infty }=
\bigcup_{n=1}^\infty S_n$. Note that $S_{\infty}$ is a
regular ring.

\begin{prop}
\label{ibn2}
Let $S_1,S_2, \dots, S_\infty$ be as above. The ring $S_{\infty }$ is
simple, regular, and purely infinite, and $K_0(S_{\infty})$ is an
infinite cyclic group, with generator $[S_\infty]$.
\end{prop}

\begin{proof} Since each $S_n$ is regular (Theorem \ref{main}), so is
$S_\infty$. Next, notice that $S_{\infty}$ is not a division ring,
because it contains at least one infinite idempotent. Let $\alpha$ be a
nonzero element of $S_\infty$. Then $\alpha \in S_n$ for some $n$, and
the proof of Proposition \ref{soc}(a) shows that there exist $s,t \in
S_n$ such that $s\alpha t= 1- \sum_{i=0}^n y_ix_i$. Consequently,
$(x_{n+1}s) \alpha (ty_{n+1}) =1$. Thus, $S_{\infty }$ is
simple and purely infinite because of Theorem \ref{pinfthm}.

We finally compute $K_0(S_{\infty})$. By Corollary \ref{KzeroS},
$K_0(S_n)\cong
\Z$ for all $n$, with $[S_n] \mapsto 1$. Proceeding as in Theorem
\ref{inflim1}, we obtain
$K_0(S_{\infty })\cong \varinjlim K_0(S_n) =\Z$, with
$[S_\infty] \mapsto 1$.
\end{proof}

Theorem \ref{K0} and Proposition \ref{ibn2} provide us with examples
of purely infinite simple regular rings whose $K_0$ groups are cyclic
of arbitrary order. Using these rings as basic building blocks, we can
construct purely infinite simple regular rings whose $K_0$'s are
arbitrary countable abelian groups -- see Section 8.

\section{The Rosenmann-Rosset and Schofield constructions}

In this section we consider two constructions associated to the free
algebra over a field $k$, and we will prove that they are isomorphic
purely infinite simple regular $k$-algebras. In Section \ref{isoQT}, we
will relate these algebras to the ones constructed in Section 5.
Throughout this section, $R$ will denote a free $k$-algebra on the
finite alphabet
$\{x_0,\dots ,x_n\}$ with $n> 0$.

First, we briefly quote an example of Rosenmann and Rosset
\cite[Section 3]{R-R}. For the general theory of rings of quotients,
we refer the reader to \cite{Stenstrom}. Let $\mathcal{F}_{fc}$ be
the Gabriel topology whose basic neighborhoods of $0$ are the right
ideals of finite codimension; see \cite[Theorem 1.1]{R-R}. Consider
the $\mathcal{F}_{fc}$-localization of $R$, denoted by $R_{fc}$. Then
$R_{fc}$ is a ring and we can associate to each right $R$-module a
right $R_{fc}$-module $M_{fc}$ by the rule $M_{fc}=\varinjlim _{I\in
\mathcal F_{fc}}\Hom (I,M)$. (Note that since every right ideal in
$R$ is projective we have $M_{fc}=M_{\mathcal F _{fc}}$ by
\cite[Lemma 2.3]{R-R}.) Recall that a module $M$ is {\it torsion}
(with respect to $\mathcal F _{fc}$) if for every $x\in M$, its
annihilator, $\text{ann}_R(x)$ is in $\mathcal F _{fc}$. If $M$ is
torsion then $M_{fc}=0$. Clearly, $M$ is torsion if and only if all
its cyclic submodules are finite dimensional.

We recall one of the main results in \cite{R-R}:

\begin{theor}
\label{RR1} {\rm \cite[Theorem 5.1]{R-R}} Let $R$ and $R_{fc}$ be as
above. Then, for $\ell,m \in \N$, we have $R_{fc}^\ell\cong
R_{fc}^m$ if and only if $\ell\equiv m \pmod{n}$.
\end{theor}

Another construction was quoted by A. H. Schofield in a private
communication \cite{SchDicks}. We thank him for allowing us to
present this construction here.

 We denote by $\mathcal{P}$ the
category of finitely generated projective right $R$-modules. Let
$\Phi =\text{Mor}(\mathcal{P})$ denote the class of all homomorphisms
between finitely generated projectives.

\begin{prop} {\rm \cite{SchDicks}} \label{sch1} Let $R$ be as above.
Let $\Sigma $ be the family of
all \underbar{mono}morphisms in $\Phi$ whose images have finite
codimension.

{\rm (a)} $\Sigma$ is composition-closed: If $f,g\in \Sigma$ and the
composition $fg$ is defined, then $fg\in \Sigma$.

{\rm (b)} Every $f\in \Sigma$ is a non-zero-divisor in $\Phi$: If
$g\in \Phi$ and $fg$ is defined (respectively, $gf$ is defined),
then $fg=0$ (respectively, $gf=0$) implies $g=0$.

{\rm (c)} $\Sigma$ satisfies the right Ore condition relative to
$\Phi$: If $f\in \Sigma$ and $g\in \Phi$ with the same codomain, then
there exist $f'\in \Sigma$ and $g'\in \Phi$ such that $gf'=fg'$ (with
appropriate conditions on domains and codomains).
\end{prop}

\begin{proof} Part (a) is clear.

(b) Let $f\in
\Sigma$ and $g\in \Phi$. Since $f$ is injective it is clear that
$fg=0$ implies $g=0$. Now assume that $gf=0$. Let $Q$ be the
codomain of $f$. Then $Q$ is a finitely generated free $R$-module,
and since
$\im (f)\subseteq \ker (g)$, we see that $\ker (g)$ has finite
codimension. Hence, $\im (g)$ is finite-dimensional and a free
$R$-module, so it is zero.

(c) Now assume that $f:X\rightarrow Y$ is in $\Sigma$ and that
$g:Z\rightarrow Y$ is in $\Phi$. Form the pullback
\begin{equation}\begin{CD}
P @>{f'}>> Z\\
@V{g'}VV @VV{g}V\\
X @>{f}>> Y
\end{CD} \notag \end{equation}
Then
$gf'=fg'$. We will show that $P \in \mathcal{P}$ and $f'\in \Sigma$,
which will finish the proof. As is well known, since $f$ is
injective it follows that
$f'$ is injective (e.g., \cite[Exercise 2.47]{Rotman}). On the other
hand
$f'(P)=g^{-1}(f(X))$. Thus
$Z/f'(P)$ embeds in
$Y/f(X)$, which is finite dimensional, and so $f'(P)$ has finite
codimension in $Z$. Since $f'$ is injective, $P$ is free of finite
rank by
\cite[Theorem 4]{Lew}, and therefore $f'\in \Sigma$.
\end{proof}

Continue with $R$ and $\Sigma$ as above. Let $\mathcal{C}$ be the set
of all pairs
$(f,s)$ with
$f\in
\Phi$ and
$s\in \Sigma$ such that $f$ and $s$ have the same domain. We define a
relation on $\mathcal{C}$. Let $(f,s)$ and $(f',s')$ be in
$\mathcal{C}$, with
\begin{equation}\begin{CD} Q @<{f}<< P @>{s}>> U \qquad\qquad
\mbox{and} \qquad\qquad Q' @<{f'}<< P' @>{s'}>> U'.
\notag\end{CD}\end{equation}
Then $(f,s)\sim (f',s')$ if and only if $Q=Q'$ and
$U=U'$ and there are maps $h:U''\rightarrow P$ and $h':U''\rightarrow
P'$ such that $sh=s'h'\in \Sigma$ and $fh=f'h'$. By \cite[Section
10.3]{Weibel}, $\sim$ is an equivalence relation, and we get a
quotient category $\mathcal{P}\Sigma^{-1}$. Now let $Q$ be the ring
consisting of all equivalence classes $[(f,s)]$ with $(f,s) \in C$ and
$f,s:P\rightarrow R$. Notice that
$Q=\mbox{End}_{\mathcal{P}\Sigma^{-1}}(R)$. There is a canonical map
$\phi :R\rightarrow Q$ defined by $\phi (a)=[(a,1)]$. The full
subcategory of Mod-$Q$ consisting of all the induced finitely
generated projective modules is equivalent to the category with the
same objects as $\mathcal{P}$ and morphisms $\mathcal{P}\Sigma^{-1}$.
The map $\phi:R\rightarrow Q$ gives the universal localization of $R$
with respect to $\Sigma$. Note that since $\Sigma $ is a
non-zero-divisor class in $\Phi$, the map $\phi$ is injective.

Since ${_RQ}$ is the direct limit of $\Hom _R(eR^m,R)s^{-1}$, where
$s:eR^m\rightarrow R$ is in $\Sigma$ and $e$ ranges through all the
idempotent matrices over $R$, and since $\Hom _R(eR^m,R)$ is
projective as a left $R$-module, it follows that ${_RQ}$ is flat.

\begin{theor} {\rm \cite{SchDicks}} \label{sch2} Let $R$ and
$\Sigma$ be as above. Then the universal
localization
$Q=R_{\Sigma}$ is a purely infinite simple regular ring.
\end{theor}

\begin{proof} Let $N$ be a finitely presented
right $R$-module. By \cite[Theorem 2]{Lew}, there is a finitely
generated free submodule $L$ of $N$ of finite codimension. There is
an exact sequence
$$0\rightarrow R^\ell\rightarrow R^m\rightarrow
N/L\rightarrow 0 ,$$
and so the map $R^\ell\rightarrow R^m$ must be
in $\Sigma$, which gives $(N/L)\otimes _RQ=0$. Now since ${_RQ}$ is
flat, we have an exact sequence
$$0\rightarrow L\otimes
_RQ\rightarrow N\otimes _RQ\rightarrow (N/L)\otimes _RQ\rightarrow
0.$$
Since $(N/L)\otimes _RQ=0$ we conclude that $N\otimes _RQ$ is a
free $Q$-module. Since every finitely presented $Q$-module is induced
from a finitely presented $R$-module, we conclude that $Q$ is a
regular ring such that every finitely generated projective module is
free. Thus $Q$ is simple and purely infinite by Example
\ref{pinfexamples}(d).
\end{proof}

\begin{theor}
\label{RRS} Let $R$, $R_{fc}$ and $Q$ be as above. Then, $R_{fc}\cong
Q$ as $k$-algebras.
\end{theor}

\begin{proof}
Let $s\in \Sigma$. Since finitely generated projective right
$R$-modules are free, we can assume without loss of generality that
$s: R^\ell\rightarrow R^m$ for some $\ell,m \in \mathbb{N}$, with
$N:=\mbox{coker}(s)$ finite-dimensional. We have therefore a short
exact sequence
\begin{equation} \begin{CD}
0 @>>> R^\ell @>s>> R^m @>>> N  @>>> 0. \notag
\end{CD}\end{equation}
Proceeding as in \cite[page 368]{R-R} and taking into account that
$N$ is a torsion module with respect to $\mathcal F _{fc}$ and thus
$N_{fc}= 0$, we get an
exact sequence
\begin{equation} \begin{CD}
0 @>>> R_{fc}^\ell @>s>> R_{fc}^m @>>> 0. \notag
\end{CD}\end{equation}
We conclude that every map in $\Sigma$ becomes invertible over
$R_{fc}$. Thus, using the universal property of $Q$ with respect to
$R$ and $\Sigma$, we conclude that there exists a unique $k$-algebra
homomorphism $\rho :Q\rightarrow R_{fc}$ which restricts to the
identity map on $R$. Since $Q$ is simple by Theorem \ref{sch2},
$\rho$ is an injective homomorphism. Now, if $I_R\leq R_R$ is a right
ideal of finite codimension, $i:I\rightarrow R$ is the natural
inclusion map, and $f:I\rightarrow R$ is any right $R$-module
homomorphism, then $\rho ([(f,i)])=[f]$, whence $\rho$ is an
isomorphism, as desired.
\end{proof}

Schofield also proved that $K_0(Q)\cong \Z/n\Z$. Because of
Theorem \ref{RRS} and Proposition \ref{vrgroup}, this implies the
result in Theorem
\ref{RR1} (\cite[Theorem 5.1]{R-R}). We will
give an alternate proof of this fact in Theorem \ref{po'ci}, by
using the techniques developed in Section 5.

\section{Isomorphism of $Q$ with $T$} \label{isoQT}

In this section we will prove that the algebra $Q$ (and so also the
algebra $R_{fc}$) constructed in Section 6 is isomorphic to a
particular instance of our construction in Section 5. For this, we
need to introduce the algebra of rational series, denoted by $\krat$,
associated to a finite alphabet $X$. This algebra
plays an important role in other parts of mathematics, particularly
formal language theory and the theory of codes, see \cite{BR}.

We fix some
notation for the rest of this section. Let $X= X_n= \{x_0,\dots,x_n\}$
be a finite alphabet, with $n>0$. The {\it division closure\/} of $k\la
X\ra$ in
$k\la\la X\ra\ra$ is called the {\it algebra of rational series\/} and
it is denoted by $\krat$. By definition, $\krat$ is the smallest
subalgebra of $k\la\la X\ra\ra$ containing $k\la X\ra$ and closed
under inversion, in the sense that if an element of the subalgebra is
invertible in $k\la\la X\ra\ra$, then it is already invertible in the
subalgebra. The algebra $\krat$ is local, with maximal ideal
consisting of the elements with constant term $0$. By \cite[Exercise
7.1.10]{free}, a square matrix over $k\la\la X\ra\ra$ (respectively
$\krat$) is invertible if and only if its image under the
augmentation map is invertible over $k$ (see also
\cite[pp.~408-409]{Cohn2}). It follows from this that
$\krat$ coincides with the {\it rational closure\/} of $k\la X\ra$ in
$k\la\la X \ra\ra$, that is, $\krat$ is the smallest subalgebra of
$k\la\la X\ra\ra$ containing $k\la X\ra$ and such that every square
matrix with coefficients in the subalgebra which is invertible over
$k\la\la X\ra\ra$ is already invertible over the subalgebra. The
right $\tau _i$-derivations $\delta _i$ on $k\langle \langle
X\rangle\rangle $ defined in Example \ref{tyux} coincide with the
left transductions associated to the monomials $x_i$. These left
transductions are defined in \cite[page 105]{free} as follows. For a
fixed monomial $x_{i_1}\cdots x_{i_r}$ of degree $r$ define the left
transduction for this monomial as the $k$-linear map $a\mapsto a^*$
of $k\langle X\rangle $ into itself which sends any monomial of the
form $bx_{i_1}\cdots x_{i_r}$ to $b$ and all other monomials to $0$.
This map extends in the obvious way to the power series algebra
$k\langle\langle X\rangle\rangle$. By \cite[page 135]{free}, the
algebra of rational series $\krat $ is invariant under all
transductions. We conclude that $\krat$ is a local subalgebra of
$k\la\la X\ra \ra$ containing $k\la X\ra$ and invariant under the
$\delta _i$. Therefore we can build the algebras $S= S_n= \krat \la
Y;\tau ,\delta \ra$ and $T= T_n= S/I=\krat _f$ based on $\krat$; see
Section 5.

An important property of $\krat $ is that it is a universal
localization of $k\la X\ra$. Although this fact seems to be well
known, we do not have an explicit reference, and so we sketch a proof
as follows.

\begin{prop}\label{kratuniv} Let $\Sigma'$ be the set of those square
matrices over $k\la X\ra$ which become invertible over $k\la\la
X\ra\ra$. Then
$\krat$ is the universal localization of $k\la X\ra$
with respect to $\Sigma'$.
\end{prop}

\begin{proof} Note that $\Sigma '$ is the
set of matrices of the form $A_0+B$, where $A_0$ is an invertible
matrix over $k$ and all the entries of $B$ have constant term $0$.
Since $\krat$ is the rational closure of $k\la X\ra$ in $k\la\la X
\ra\ra$, we have a unique algebra homomorphism $\phi: k\la X\ra
_{\Sigma ' }\rightarrow \krat$ which is the identity on $k\la X\ra$.
It is easy to see from \cite[Theorem 7.1.2]{free} that $\phi$ is
surjective. (This holds in general when we consider the rational
closure of an inclusion of rings.)

By \cite[Lemma 5.9.4]{free}, the inclusion map $k\la X\ra \rightarrow
k\la\la X\ra\ra$ is honest (i.e., it sends full matrices to full
matrices), and hence so is the inclusion map
$k\la X\ra \rightarrow \krat$. Since $\Sigma'$ is a factor-closed,
multiplicative set of matrices over $k\la X\ra$, it now follows from
\cite[Proposition 7.5.7(ii)]{free} that $\phi$ is injective. Therefore
$\phi$ is an isomorphism.
\end{proof}

By a result of Bergman and Dicks
(\cite[Theorem 4.9]{scho}), every universal localization of a
hereditary ring is also hereditary. In particular, $\krat$ is a
hereditary ring, and indeed it is a fir by \cite[Theorem
7.10.7]{free}. The algebra $k\la\la X\ra\ra$ is just a semifir, but
not a fir; see \cite[Proposition 2.9.18 and Proposition
5.10.9]{free}. We summarize the results obtained in Section 5 for the
algebras $S$ and $T$ based on the algebra of rational series $\krat$,
along with the facts just mentioned.

The following notation will be
helpful. If $a_0,\dots,a_n$ are elements in a $k$-algebra $V$, let
$\Sigma'_k(a_0,\dots,a_n)$ denote the set of square matrices over $V$
of the form $A_0+B$ where $A_0$ is an invertible matrix over $k$ and
the entries of $B$ consist of noncommutative polynomials with zero
constant term evaluated at $(a_0,\dots,a_n)$. (Note that to check
whether all such matrices are invertible over $V$, it suffices to
consider those for which $A_0$ is an identity matrix.) In particular,
if
$V= k\la X\ra$ and the $a_i=x_i$, then $\Sigma'_k(x_0,\dots,x_n)$
equals the set $\Sigma'$ in Proposition \ref{kratuniv}.

\begin{theor}
\label{rat} Let $\Sigma _1= \Sigma'_k(x_0,\dots,x_n) \cup \{(x_0,\dots
,x_n)^T\}$. Set $S=\krat\la Y;\tau,\delta\ra$, let $I$ be
the ideal generated by the idempotent $e=1-\sum _{i=0}^n y_ix_i$, and
write $T=S/I$. Then $T$ is the universal localization of $k\la X\ra$
with respect to $\Sigma _1$. In particular, $T$ is a hereditary
algebra.

To write the universal property of $T$ in more elementary form,
suppose that $V$ is a $k$-algebra containing elements
$a_0,\dots,a_n$ and
$b_0,\dots,b_n$ such that
\begin{enumerate}
\item [(1)] $a_ib_j= \delta_{ij}$ for all $i,j$.
\item [(2)] $b_0a_0+ b_1a_1+ \dots+ b_na_n =1$.
\item [(3)] All matrices in $\Sigma'_k(a_0,\dots,a_n)$ are invertible
over $V$.
\end{enumerate}
Then there exists a unique $k$-algebra homomorphism $\psi:
T\rightarrow V$ such that $\psi(x_i)=a_i$ for all $i$ and
$\psi(y_j)=b_j$ for all $j$.
\end{theor}

\begin{proof} That $T= k\la X\ra_{\Sigma_1}$ follows from
Proposition \ref{kratuniv} and
Theorem \ref{main}. The fact that $T$ is hereditary follows from
\cite[Theorem 4.9]{scho}.
\end{proof}

Now we bring the construction of Theorem \ref{sch2} into play, but
with a slight change of notation. Let $Y=\{y_0,\dots ,y_n\}$ and
set $R= k\langle Y\rangle$. We may identify $R$ with the
$k$-subalgebra of $S$ generated by $Y$. If $r$ is any nonzero
element of $R$ and $y_J$ is a monomial of maximum length occurring
in $r$, then $x_{J^*}r$ is an element of $k\langle X\rangle$ with
nonzero constant term. Hence, $x_{J^*}r$ is invertible in $S$, and
so $r\notin I$. Therefore $R\cap I= 0$, and thus we can identify $R$
with its image in $T$ under the quotient map.

Let $Q$ be the universal localization of
$R$ with respect to the set $\Sigma$ of monomorphisms between
finitely generated projective right $R$-modules with
finite dimensional cokernels. We want to
prove that all maps in $\Sigma$ are invertible over $T$.

We need the following fact:

\begin{lem}
$T$ is flat as a left $R$-module.
\end{lem}

\begin{proof}
We first check that $_RS$ is free. To see this, let $\{p_\alpha\}$ be
a $k$-basis of
$\krat$. Then the elements of
$S$ can be uniquely written as $\sum \lambda_{I,\alpha} y_Ip_\alpha$.
We have $$S= \bigoplus_I y_I \krat= \bigoplus_{I,\alpha}
ky_Ip_\alpha= \bigoplus_\alpha \bigl( \bigoplus_I ky_I \bigr)
p_\alpha= \bigoplus_\alpha Rp_\alpha,$$ so $_RS$ is free. Since $S$
is regular, $_ST$ is flat, hence $_RT$ is flat.
\end{proof}

\begin{prop}
\label{sigmainv} Every map in $\Sigma$ is invertible over $T$, and so
there is a unique $k$-algebra homomorphism $Q\rightarrow T$ which is
the identity on $R$.
\end{prop}

\begin{proof} Notice that every finitely generated projective
right $R$-module is free. Consider a short exact sequence
\begin{equation} \begin{CD}
0\rightarrow R^i @>g>> R^j\longrightarrow N\longrightarrow 0 \end{CD}
\notag
\end{equation}
with $g\in \Sigma$. Since $_RT$ is flat, we have an exact sequence
\begin{equation} \begin{CD}
0\rightarrow T^i @>{g}>> T^j \longrightarrow N\otimes_R T
\longrightarrow 0,\end{CD} \tag{*}
\end{equation}
and it is enough to see that $N\otimes_RT =0$ to get the desired
result.

Since $g\in\Sigma$ we have  $\dim_k(N) <\infty$, and so $J:=
\text{ann}_R(N)$ is an ideal of $R$ of finite codimension. Since $T$
is regular, the sequence (*) splits, so $N\otimes_R T$ is a
projective right $T$-module; in particular, it embeds in $T_T$. For
$\gamma \in N$, $(\gamma \otimes 1)R$ has finite $k$-dimension since
$\gamma J=0$ and $\dim_k(R/J) <\infty$. So it is enough to show that
$T$ does not contain nonzero finite-dimensional right $R$-submodules.
We will prove this in the following form: If $\alpha \in S\setminus
I$, then
$\alpha R$ contains an infinite family of elements which are linearly
independent modulo $I$.

If $\alpha \in S\setminus I$, then, as shown in the
proof of Theorem \ref{main}, we can find
$m\in X^*$ such that $m\alpha \in \krat$ and $m\alpha \ne 0$. By Lemma
5.1, there is a word $w\in Y^*$ such that the product $p= m\alpha w$
is an invertible element of $\krat$. As noted
earlier in the section, $R\cap I= \{0\}$, whence
$1,y_1,y_1^2,\dots$ are $k$-linearly independent modulo $I$. Thus, the
sequence $p^{-1}m\alpha w, p^{-1}m\alpha wy_1, p^{-1}m\alpha wy_1^2
,\dots$ is linearly independent modulo $I$. It follows immediately that
$\alpha w,
\alpha wy_1,
\alpha wy_1^2 ,\dots$ is linearly independent modulo $I$, as desired.

Therefore all maps in $\Sigma$ are indeed invertible over $T$, so we
get a $k$-algebra homomorphism $\psi: Q \rightarrow  T$ which is the
identity on $R$.
\end{proof}

We are now ready to prove the main result of this section.

\begin{theor}
\label{QisoT} Let $T$ be the universal localization of $k\la X \ra$
with respect to the set $\Sigma _1$ described in Theorem \ref{rat}.
Let $Q$ be the universal localization of $R=k\la Y\ra$ with respect
to the set $\Sigma$ described in Section 6. Then there is a unique
$k$-algebra isomorphism $\psi :Q\rightarrow T$ which is the identity
on $R$.
\end{theor}

\begin{proof}
By Proposition \ref{sigmainv}, there exists a unique $k$-algebra
homomorphism $\psi :Q\rightarrow T$ which is the identity on $R$.
Since $Q$ is simple, the map $\psi $ must be injective. Now $\psi$
sends the row matrix $(y_0,\dots ,y_n)$ over $R$ to the row matrix
$(y_0,\dots ,y_n)$ over $T$, which is the inverse of the column
$(x_0,\dots, x_n)^T$ over $T$. The row $(y_0,\dots ,y_n)$ is invertible over $Q$,
so there exists a column $(x_0',\dots ,x_n')^T$ over $Q$ which is
the inverse of $(y_0,\dots ,y_n)$, and obviously
$\psi (x_i')=x_i$ for all $i$. Let $A$ be an $m\times m$ matrix in
$\Sigma'_k(x_0',\dots ,x_n')$. Then $\psi(A)$ is invertible in
$M_m(T)$, and so, $\psi$ being injective,
$A$ must be a non-zero-divisor in $M_m(Q)$. Since $Q$ is regular by
Theorem
\ref{sch2}, we get that $A$ is invertible in $M_m(Q)$. By the universal
property of
$T=k\la X\ra _{\Sigma _1}$, there exists a unique $k$-algebra
homomorphism $\varphi :T\rightarrow Q$ sending $(x_0,\dots ,x_n)^T$
to $(x_0',\dots ,x_n')^T$. The row matrix $(y_0,\dots ,y_n)$, being
the inverse of $(x_0,\dots ,x_n)^T$ in $T$ as well as the inverse of
$(x_0',\dots ,x_n')^T$ in $Q$, must be sent to itself by $\varphi$.
Therefore we conclude that
$\psi $ and
$\varphi$ are mutually inverse isomorphisms.
\end{proof}

As a consequence of our results in Section 5, we can derive an
alternate proof of Schofield's result that $K_0(Q)\cong \Z/n\Z$.
As we observed at the end of Section 6, this implies in turn the
Rosenmann-Rosset result \cite[Theorem 5.1]{R-R}.

\begin{theor} \label{po'ci}
Let $Q$ be as above. Then $K_0(Q)\cong \Z/n\Z$, with $[Q]
\mapsto [1]_n$.
\end{theor}

\begin{proof}
This follows immediately from Theorems \ref{QisoT} and
\ref{K0}.
\end{proof}

\section{Realizing groups as $K_0$ of purely infinite simple regular
rings}

R\o rdam has proved that all countable abelian groups appear as
$K_0$'s of purely infinite simple C*-algebras \cite[Theorem
8.1]{rordamclassif}. In this section we prove a similar result for
purely infinite simple regular rings.  Our construction follows the
same pattern as R\o rdam's but requires much more care with the
technical details, since we have to work with more complicated
universal properties.

For a field $K$ and an integer $n\ge 2$, we will write
$S_{n,K}$ and $T_{n,K}$ for the $K$-algebra versions of the algebras
constructed in Section
\ref{isoQT}. Namely,
$$S_{n,K}=K_{\mbox{rat}}\langle X_n\rangle\langle Y_n;\tau,
\delta \rangle \, ,$$
where $X_n=\{x_0,x_1,\dots ,x_n\}$ and
$Y_n=\{y_0,y_1,\dots ,y_n\}$, and $T_{n,K}=S_{n,K}/I_{n,K}$,
where $I_{n,K}$ is the ideal of $S_{n,K}$ generated by the idempotent
$e_n :=1-\sum _{p=0}^n y_px_p$. By Theorems \ref{main} and \ref{K0},
$T_{n,K}$ is a purely infinite simple regular ring, and $K_0(T_{n,K})
\cong
\Z/n\Z$ with $[T_{n,K}] \mapsto [1]_n$. We denote by
$T_{0,K}$ the inductive limit of the sequence $(S_{n,K})_{n\ge 2}$
along the natural inclusions $S_{n,K} \hookrightarrow S_{n+1,K}$. Note
that $T_{0,K}$ contains the idempotents $e_0,e_1,\dots$ and that
$$ x_ie_n = \begin{cases} 0 &\quad (i\le n)\\ x_i
&\quad (i>n)\, , \end{cases} \qquad\qquad\qquad
e_ny_j = \begin{cases} 0 &\quad (j\le n)\\ y_j &\quad (j>n)\, .
\end{cases}$$
By Proposition \ref{ibn2}, $T_{0,K}$ is simple, regular, and
purely infinite. Moreover, $K_0(T_{0,K}) \cong \Z$ with $[T_{0,K}]
\mapsto 1$.

We shall require the universal property for $T_{n,K}$, $(n\ge 2)$,
given in Theorem \ref{rat}. Analogously,
$T_{0,K}$ possesses the following universal property:

\begin{lem}\label{T0Kuniv} Let $V$ be a $K$-algebra containing elements
$a_0,a_1,\dots$ and $b_0,b_1,\dots$ such that $a_ib_j= \delta_{i,j}$
for all $i,j$, and such that all the matrices in
$\Sigma'_K(a_0,a_1,\dots)$ are invertible over
$V$. Then there exists a unique $K$-algebra homomorphism $\psi: T_{0,K}
\rightarrow V$ such that $\psi(x_i)= a_i$ for all $i$ and $\psi(y_j)=
b_j$ for all
$j$.
\end{lem}

\begin{proof} This follows from Lemmas \ref{kratuniv} and \ref{soc}(b).
\end{proof}

 We will take as canonical representatives of the nonzero cyclic
groups the groups
$\Z_m$, the integers mod $m$, where $m$ is either $0$ or an integer
larger than $1$. For uniformity of notation, we write $[a]_0= a$ for
$a\in\Z$. Every group homomorphism
$\varphi :\Z _n\rightarrow
\Z _m$ is given by multiplication by some integer $\ell$, and we
choose $\ell$ such that $1\le \ell\le m $ in case $m\ge 2$. Recall
that for any ring $R$ and $\ell>0$, there is a canonical group
isomorphism $K_0(M_\ell(R)) \cong K_0(R)$ with $[M_\ell(R)] \mapsto
\ell[R]$. Hence, we identify
$$\bigl( K_0(M_\ell(T_{m,K})),[M_\ell(T_{m,K})] \bigr) =
(\Z_m,[\ell]_m)$$
for all $\ell>0$ and $m\in \{0\} \cup \{2,3,\dots\}$.

It will be convenient to set $e_{-\ell}= e_\ell= 1-y_0x_0-y_1x_1-\cdots
-y_\ell x_\ell$ for all $\ell\ge 0$. If $\ell\le 0$, then by
$M_\ell(T_{0,K})$ we will understand the corner ring $e_\ell T_{0,K}
e_\ell$. In this case
$K_0(M_\ell(T_{0,K})) \cong K_0(T_{0,K})$ with $[M_\ell(T_{0,K})]
\mapsto [e_\ell T_{0,K}]$. Note that $y_0x_0, \dots,
y_{-\ell}x_{-\ell}$ are pairwise orthogonal idempotents equivalent to
1, and also orthogonal to $e_\ell$. Since $e_\ell+ y_0x_0+ y_1x_1+
\dots+ y_{-\ell}x_{-\ell} =1$, we have
$[e_\ell T_{0,K}] +(-\ell+1)[T_{0,K}]= [T_{0,K}]$ in $K_0(T_{0,K})$,
whence
$[e_\ell T_{0,K}]= \ell[T_{0,K}]$. Hence, we can make the
identification
$$\bigl( K_0(M_\ell(T_{0,K})),[M_\ell(T_{0,K})] \bigr) =
(\Z_0,[\ell]_0),$$
in parallel with the previous identifications.

The following lemma and corollary give the key for our construction.

\begin{lem}
\label{lemreal} Let $\varphi :\mathbb Z_n\rightarrow \mathbb Z_m$ be a
group homomorphism given by multiplication by $\ell$, where we adopt
the above conventions. Let $K$ be a field and $t$ an indeterminate.
Then there exists a
$K$-algebra homomorphism $\psi :T_{n,K}\rightarrow M_\ell(T_{m,K(t)})$
such that $K_0(\psi )=\varphi$.
\end{lem}

\begin{proof} With the canonical identifications we have made for the
$K_0$ groups, it is automatic that $K_0(\psi)=\varphi$ will hold for
any  $K$-algebra homomorphism $\psi :T_{n,K}\rightarrow
M_\ell(T_{m,K(t)})$. Thus, only the existence of such homomorphisms
needs to be established.
We will consider several cases, in the first two of which $\ell>0$.

{\bf Case 1}. Assume that $n,m\ge 2$.

We can write $\ell n= hm$ for some positive integer
$h$. There are elements $x_0',\dots, x_{hm}'$ and
$y_0',\dots ,y_{hm}'$ in
$T_{m,K(t)}$ implementing an isomorphism $T_{m,K(t)}\cong
T_{m,K(t)}^{hm +1}$ and with each $x_i'$ being a nontrivial product
of the standard elements $x_i$ of $T_{m,K(t)}$. For example, take
$x_0',\dots, x_{hm}'$ and
$y_0',\dots ,y_{hm}'$ to be the sequences
\begin{align}
x_0^h \, , &x_1x_0^{h-1} \, , \dots \, , x_mx_0^{h-1} \, ,
x_1x_0^{h-2} \, , \dots \, , x_mx_0^{h-2} \, , \dots \, , x_1 \,
,\dots \, ,x_m \qquad \text{and} \notag\\
y_0^h \, , &y_0^{h-1}y_1 \, , \dots \, , y_0^{h-1}y_m \, , y_0^{h-2}y_1
\, , \dots \, , y_0^{h-2}y_m \, , \dots \, , y_1 \, ,\dots \, ,y_m \,
, \notag
\end{align}
respectively. Define matrices
$A_i,B_j\in M_\ell(T_{m,K(t)})$ for $i,j=0,1,\dots ,n$ as follows:
\begin{gather} A_0=\mbox{diag}(t,\dots,t,x_0')\, ,\qquad\qquad\qquad
(A_i)_{\alpha ,\beta}=\delta _{\beta,\ell}x'_{(i-1)\ell+\alpha} \quad
(1\le i\le n)\, ,\notag\\
B_0=\mbox{diag}(t^{-1},\dots ,t^{-1},y_0')
\, ,\qquad\qquad (B_j)_{\alpha ,\beta }=\delta
_{\alpha,\ell}y'_{(i-1)\ell+\beta}\quad (1\le j\le n)\, .
\notag\end{gather}
 It is easy to check that $A_iB_j= \delta_{i,j}I$ for all $i,j$ and
$\sum_{i=0}^n B_iA_i= I$.

Observe that if $p(Z_0,\dots,Z_n)$ is a
noncommutative polynomial over $K$ with zero constant term, then
$p(A_0,\dots,A_n)$ is a matrix in $M_\ell(T_{m,K(t)})$ whose entries
come from either $tK[t]$ or $K[t]\la X_m\ra X_m$ (where the latter
notation refers to the left ideal of $K[t]\la X_m\ra$ generated by
$X_m$). Consequently, any $r\times r$ matrix in
$\Sigma'_K(A_0,\dots,A_n)$, when viewed as a block form of an $\ell
r\times \ell r$ matrix over $T_{m,K(t)}$, consists of a sum
$C_0+C_1+C_2$ where $C_0$ is an invertible matrix over $K$, all
entries of $C_1$ lie in $tK[t]$, and all entries of $C_2$ lie in
$K[t]\la X_m\ra X_m$. Now $C_0+C_1$ is an invertible matrix over
$K(t)$ (note that we
need the variable $t$ to achieve this statement), whence $C_0+C_1+C_2$
lies in $\Sigma'_{K(t)}(x_0,\dots,x_m)$. It follows that every matrix
in $\Sigma'_K(A_0,\dots,A_n)$ is invertible over
$M_\ell(T_{m,K(t)})$. By the universal property of the
algebras $T_{n,K}$ (Theorem \ref{rat}), there is a unique
$K$-algebra homomorphism $T_{n,K}\rightarrow M_\ell(T_{m,K(t)})$
sending $(x_i)$ to $(A_i)$ and $(y_j)$ to $(B_j)$.

{\bf Case 2}. Assume that $n=0$ and $\ell>0$ (with $m$ arbitrary).

This case is
easier than the previous one. Just take an infinite sequence of
diagonal matrices
$A_i$, $B_j$ in
$M_\ell(T_{m,K})$, where the $A_i$'s consist of nontrivial products of
the $x$'s in the diagonal, the $B_j$'s consist of nontrivial products
of the $y$'s in the diagonal, and they satisfy the rules
$A_iB_j=\delta _{ij}$. For example, take
$$A_i= \mbox{diag} (x_0x_1^i,\dots,x_0x_1^i) \qquad\quad \text{and}
\qquad\quad B_j= \mbox{diag} (y_1^jy_0,\dots,y_1^jy_0).$$
for $i,j= 0,1,\dots$. Then we can identify $\Sigma'_K(A_0, A_1,\dots)$
with a subset of $\Sigma'_K(x_0,x_1)$, and so every matrix in
$\Sigma'_K(A_0, A_1,\dots)$ is
invertible over $M_\ell(T_{m,K})$. By the universal property of
$T_{0,K}$ (Lemma \ref{T0Kuniv}), there is a  $K$-algebra homomorphism
$\psi:T_{0,K}\rightarrow M_\ell(T_{m,K})$. (Here the use of the new
variable $t$ is not necessary.)

{\bf Case 3}. Now assume that $n=0$ and $\ell\le 0$. Necessarily,
$m=0$.

Recall that in this case our convention is to set
$M_\ell(T_{0,K})= e_\ell T_{0,K}e_\ell$. Define the following elements
$A_i$, $B_j$ in
$M_\ell(T_{0,K})$ for $i,j\ge 0$:
$$A_i= e_\ell x_{-\ell+1+i} \qquad\qquad \mbox{and} \qquad\qquad B_j=
y_{-\ell+1+j} e_\ell.$$  Then one has
$A_iB_j=\delta _{ij} e_\ell$ for all $i,j$. If $C$ is a
matrix in $\Sigma'_K(A_0,A_1,\dots)$, then $C= \etil_\ell\Gamma=
\etil_\ell\Gamma\etil_\ell$ for some matrix
$\Gamma$ in
$\Sigma'_K(x_{-\ell+1}, x_{-\ell+2},
\dots)$, where $\etil_\ell= \mbox{diag}(e_\ell, \dots, e_\ell)$. Note
that $\Gamma$ is invertible over
$K_{\mbox{rat}}\langle x_{-\ell+1}, x_{-\ell+2}, \dots\rangle$, with
$\etil_\ell \Gamma^{-1}= \etil_\ell \Gamma^{-1} \etil_\ell$, and so $C$
is invertible over
$M_\ell(T_{0,K})$ with inverse $\etil_\ell \Gamma^{-1}$. Thus in this
case too we get our desired
 homomorphism
$T_{0,K}\rightarrow M_\ell(T_{0,K})$.

{\bf Case 4}. Finally, suppose that $m=0$ and $n\ge 2$.
Then necessarily $\ell=0$.

Set $E= \mbox{diag}(e_n,1,1,\dots,1)$ in
$M_{n+1}(T_{0,K(t)})$, and observe that $E \sim e_0$. Hence,
$M_0(T_{0,K(t)})$ is isomorphic to the $K(t)$-algebra $M= E
M_{n+1}(T_{0,K(t)}) E$, and so it suffices to produce a $K$-algebra
homomorphism $T_{n,K} \rightarrow M$. Let us index rows and columns of
matrices in $M_{n+1}(T_{0,K(t)})$ by $0,1,\dots,n$.

Now consider matrices $A_i$, $B_j$ in $M$ for $i,j= 0,1,\dots,n$ where
\begin{gather}
A_0= \mbox{diag}(te_n,x_0,\dots,x_0)\, ,\qquad\qquad
(A_i)_{\alpha,\beta}= \delta_{\beta,i}\begin{cases} te_n &\quad
(\alpha=0)\\ x_\alpha &\quad (\alpha\ge1) \end{cases} \quad (1\le i\le
n)\, ,\notag\\
B_0= \mbox{diag}(t^{-1}e_n,y_0,\dots,y_0)\, ,\qquad\quad
(B_j)_{\alpha,\beta}= \delta_{\alpha,j}\begin{cases} t^{-1}e_n &\quad
(\beta=0)\\ y_\beta &\quad (\beta\ge1) \end{cases} \quad (1\le j\le
n)\, . \notag
\end{gather}
Then we have $A_iB_j= \delta_{i,j}E$ for all $i,j$ and
$\sum_{i=0}^n B_iA_i= E$. If $p(Z_0,\dots,Z_n)$ is a
noncommutative polynomial over $K$ with zero constant term, then
$p(A_0,\dots,A_n)= EC= ECE$ for some matrix $C\in
M_{n+1}(T_{0,K(t)})$ whose entries come from either $tK[t]\la X_n\ra$
or $K\la X_n\ra X_n$, and where $C_{i0}=0$ for $i=1,\dots,n$.
Consequently, any matrix $D$ in
$\Sigma'_K(A_0,\dots,A_n)$ can be written as $D= \Etil\Delta=
\Etil\Delta\Etil$ for some
$\Delta$ in $\Sigma'_{K(t)}(x_0,\dots,x_n)$, where $\Etil=
\mbox{diag}(E,\dots,E)$. Now
$\Delta$ is invertible over
$K(t)_{\mbox{rat}}\la X_n\ra$, and we compute that $\Etil\Delta^{-1}=
\Etil\Delta^{-1}\Etil$. Thus $D$ is invertible over $M$, with inverse
$\Etil\Delta^{-1}$. This shows that every matrix in
$\Sigma'_K(A_0,\dots,A_n)$ is invertible over $M$. Therefore by the
universal property of $T_{n,K}$, there is a
 $K$-algebra homomorphism
$T_{n,K} \rightarrow M$, and we are done.
\end{proof}

\begin{corol} \label{correal} Let $K$ be a field and $t$ an
indeterminate. Let $R$ be a $K$-algebra Morita equivalent to some
$T_{n,K}$ and $S$ a $K(t)$-algebra Morita equivalent to some
$T_{m,K(t)}$, where $n,m\in \{0\} \cup \{2,3,\dots\}$. Let $\varphi:
K_0(R) \rightarrow K_0(S)$ be a group homomorphism such that
$\varphi([R])= [S]$. Then there exists a $K$-algebra homomorphism
$\psi: R\rightarrow S$ such that $K_0(\psi)= \varphi$. \end{corol}

\begin{proof} There exists a finitely generated projective right
$T_{n,K}$-module $A$ such that $R \cong \End(A)$, and $[A]=
\ell[T_{n,K}]$ for some $\ell\in\Z$, where we may assume $1\le \ell\le
n$ in case $n\ge 2$. If $\ell>0$, then $A\cong T^\ell_{n,K}$ by
Proposition \ref{vrgroup}, whence $R\cong M_\ell(T_{n,K})$. If $\ell\le
0$, then $n=0$ and $[A]= [e_\ell T_{0,K}]$, in which case $R\cong
e_\ell T_{0,K} e_\ell = M_\ell(T_{0,K})$. Hence, there is no loss of
generality in assuming that $R= M_\ell(T_{n,K})$. Likewise, we may
assume that $S= M_{\ell'}(T_{m,K(t)})$ for some $\ell' \in\Z$, where
$1\le \ell'\le m$ in case $m\ge 2$.

As above, we can make the identifications
\begin{gather}
\bigl( K_0(T_{n,K}), [T_{n,K}] \bigr)= (\Z_n,[1]_n) \notag\\
(K_0(R),[R])= \bigl(
K_0(T_{n,K}),\ell[T_{n,K}] \bigr)= (\Z_n,[\ell]_n) \notag\\
\bigl( K_0(T_{m,K(t)}), [T_{m,K(t)}] \bigr)= (\Z_m,[1]_m) \notag\\
(K_0(S),[S])= \bigl( K_0(T_{m,K(t)}),\ell'[T_{m,K(t)}] \bigr)=
(\Z_m,[\ell']_m). \notag
\end{gather}
Then $\varphi$ is identified with a homomorphism $\Z_n
\rightarrow
\Z_m$ given by multiplication by an integer $h$, where we can assume
that $1\le h\le m$ in case $m\ge 2$, and $[h\ell]_m= [\ell']_m$.
Finally, identify
$$\bigl( K_0(M_h(T_{m,K(t)})), [M_h(T_{m,K(t)})] \bigr)= \bigl(
K_0(T_{m,K(t)}), h[T_{m,K(t)}] \bigr)= (\Z_m, [h]_m).$$

By Lemma \ref{lemreal}, there exists a $K$-algebra homomorphism
$\psi_0 : T_{n,K} \rightarrow M_h( T_{m,K(t)})$ such that
$K_0(\psi_0)= \varphi$. If $\ell>0$, then $\psi_0$ induces a
$K$-algebra homomorphism $\psi_1: R\rightarrow M_\ell( M_h(
T_{m,K(t)}))$ such that $K_0(\psi_1)= \varphi$. There is a nonzero
finitely generated projective right $S$-module $B$ such that $[B]=
[h]_m$ in $K_0(S)$ and $\End_S(B) \cong M_h(T_{m,K(t)})$. Since
$\ell[B]= [\ell h]_m= [\ell']_m= [S]$, it follows from Proposition
\ref{vrgroup} that $\ell B\cong S_S$, whence $M_\ell( M_h(T_{m,K(t)}))
\cong S$. On composing the latter isomorphism with $\psi_1$, we obtain
a $K$-algebra homomorphism $\psi: R\rightarrow S$ such that
$K_0(\psi)= \varphi$.

Now suppose that $\ell\le 0$, so that $n=0$ and $R= e_\ell T_{0,K}
e_\ell$. In this case, $\psi_0$ restricts to a $K$-algebra homomorphism
$\psi_1 : R\rightarrow fM_h(T_{m,K(t)})f$, where $f= \psi_0(e_\ell)$.
Since $[e_\ell T_{0,K}]= \ell [T_{0,K}]$ in $K_0(T_{0,K})$, we have
$$[fM_h(T_{m,K(t)})]= \ell[M_h(T_{m,K(t)})]= \ell[h]_m= [\ell']_m$$
in $K_0(M_h(T_{m,K(t)}))$. There is a nonzero finitely generated
projective right $M_h(T_{m,K(t)})$-module $C$ such that $[C]=
[\ell']_m$ in $K_0(M_h(T_{m,K(t)}))$ and $\End(C) \cong S$. By
Proposition 2.2, $fM_h(T_{m,K(t)}) \cong C$, whence $fM_h(T_{m,K(t)})f
\cong S$. As in the previous paragraph, on composing the latter
isomorphism with $\psi_1$, we obtain the desired $K$-algebra
homomorphism $\psi: R\rightarrow S$. \end{proof}

We are now ready to prove our realization result.

\begin{theor}
\label{thereal} Let $G$ be a countable abelian group, $u\in G$, and
$k$ any field. Then there exists a purely infinite simple regular
$k$-algebra $R$ such that
$K_0(R)\cong G$ with $[R] \mapsto u$.
\end{theor}

\begin{proof}
We can write $(G,u)$ as the inductive limit of a sequence
\begin{equation}\begin{CD}
(G_0,u_0) @>{\varphi_0}>> (G_1,u_1)
@>{\varphi_1}>> \cdots \, ,
\end{CD}\notag\end{equation}
 where each $G_n$ is a nonzero finitely
generated abelian group,
$u_n\in G_n$, and $\varphi_n: G_n \rightarrow G_{n+1}$ is a group
homomorphism such that $\varphi_n(u_n)= u_{n+1}$. For each $n$, we may
assume that
$$(G_n,u_n)= (G_{n,1},u_{n,1}) \times \cdots\times (G_{n,r(n)},
u_{n,r(n)})$$
where each $G_{n,i}$ is nonzero and cyclic.

Let $t_1,t_2,\dots$ be independent indeterminates over $k$, and set
$K_n= k(t_1,\dots,t_n)$ for $n=0,1,\dots$ (thus $K_0=k$). For each
$n$, let $\mathcal{M}_n$ denote the class of $K_n$-algebras
Morita equivalent to ones of the form $T_{*,K_n}$. Choose $R_{n,1},
\dots, R_{n,r(n)}$ in $\mathcal{M}_n$ together with identifications
$(K_0(R_{n,i}), [R_{n,i}]) =(G_{n,i},u_{n,i})$ for all $i$. Set $R_n=
R_{n,1} \times\cdots\times R_{n,r(n)}$, so that $(K_0(R_n),[R_n])=
(G_n,u_n)$. We shall construct $K_n$-algebra homomorphisms $\psi_n
:R_n \rightarrow R_{n+1}$ such that $K_0(\psi_n)= \varphi_n$ and each
of the component maps $\psi_{n,j}: R_n \rightarrow R_{n+1,j}$ is an
embedding. Then the inductive limit of the sequence
\begin{equation}\begin{CD}
R_0 @>{\psi_0}>> R_1 @>{\psi_1}>> \cdots
\end{CD}\notag\end{equation}
will be a regular $k$-algebra with $(K_0(R),[R]) \cong (G,u)$. Because
the $\psi_{n,j}$ are embeddings, we see that for any nonzero element
$a\in R_n$, there exist elements $x,y \in R_{n+1}$ such that
$x\psi_n(a)y =1$. Therefore $R$ will be a purely infinite simple ring.

It only remains to build the homomorphisms $\psi_n$. Fix $n$, write
$\varphi_n$ as a matrix of group homomorphisms
$$\eta_{i,j}: G_{n,i} \rightarrow G_{n+1,j} \, ,$$
and note that $\eta_{1,j}(u_{n,1}) +\dots+ \eta_{r(n),j}(u_{n,r(n)})
=u_{n+1,j}$ for all $j$. Hence, there exist finitely generated
projective right $R_{n+1,j}$-modules $P_1,\dots,P_{r(n)}$ such that
$[P_i]= \eta_{i,j}(u_{n,i})$ for all $i$ and $P_1\oplus \cdots\oplus
P_{r(n)} \cong R_{n+1,j}$. Since $R_{n+1,j}$ is purely infinite
simple, there exists a nonzero finitely generated projective right
$R_{n+1,j}$-module $P$ such that $[P]$ is the identity element of the
group $\mathcal{V}(R_{n+1,j})^*$ (see Proposition \ref{vrgroup}). We
can replace each $P_i$ by $P\oplus P_i$, and so we may assume that all
$P_i \ne 0$. Consequently, there exist nonzero pairwise orthogonal
idempotents $f_{1,j},
\dots, f_{r(n),j}$ in $R_{n+1,j}$ such that $f_{1,j} +\dots+
f_{r(n),j} =1$ and $[f_{i,j}R_{n+1,j}]= \eta_{i,j}(u_{n,i})$ for all
$i$.

Each corner $f_{i,j} R_{n+1,j} f_{i,j}$
belongs to $\mathcal{M}_{n+1}$, and we can make the identification
$$\bigl( K_0( f_{i,j} R_{n+1,j} f_{i,j}), [f_{i,j} R_{n+1,j} f_{i,j}]
\bigr)= \bigl( G_{n+1,j}, \eta_{i,j}(u_{n,i}) \bigr).$$
By Corollary \ref{correal}, there exist $K_n$-algebra homomorphisms
$\theta_{i,j}: R_{n,i} \rightarrow f_{i,j} R_{n+1,j} f_{i,j}$ such that
$K_0(\theta_{i,j})= \eta_{i,j}$. Since the $R_{n,i}$ are simple
algebras, the $\theta_{i,j}$ are embeddings. Consequently, the rule
$\psi_{n,j}(a)= \theta_{1,j}(a) +\dots+ \theta_{r(n),j}(a)$
defines a $K_n$-algebra embedding $\psi_{n,j}: R_n \rightarrow
R_{n+1,j}$ such that $K_0(\psi_{n,j})= (\eta_{1,j}, \dots,
\eta_{r(n),j})$. These $\psi_{n,j}$, finally, are the components for
the desired $K_n$-algebra homomorphism $\psi_n :R_n \rightarrow
R_{n+1}$.
\end{proof}

\begin{rema}
Note that we get a countable $k$-algebra $R$ in Theorem \ref{thereal}
in case we start with a countable field $k$.
\end{rema}

\section*{Acknowledgements}

Part of this work was done during a visit of the first author to the
Department of Mathematics of the University of California at Santa
Barbara, and visits of the second and third authors to the Centre de
Recerca Matem\`atica, Institut d'Estudis Catalans in Barcelona. The
three authors are very grateful to the host centers for their warm
hospitality. We also thank Mikael R\o rdam for interesting discussions
on the topic of Section 8.

\end{document}